\documentclass[a4paper,11pt]{article}
\usepackage{amsfonts, amsmath, amssymb, amscd, latexsym, graphics}

\newtheorem{thm}{Theorem}[section]
\newtheorem{prop}[thm]{Proposition}
\newtheorem{lem}[thm]{Lemma}
\newtheorem{cor}[thm]{Corollary}

\newtheorem{defn}[thm]{Definition}

\newcommand{\demo}{ {\it   Proof. }}
\def\qed{\hfill\rlap{$\sqcup$}$\sqcap$\par}
\newcommand{\Z}{\mathbb{Z}}

\title{\Large On uniform conjugators in torsion-free hyperbolic groups}

\author{O.\ Bogopolski \\ \small{Institute of Mathematics of} \\ \small{Siberian Branch of Russian Academy
of Sciences,} \\ {\small Novosibirsk, Russia}\\ {\small and D\"{u}sseldorf University, Germany} \\ \small{e-mail:
Oleg$\_$Bogopolski@yahoo.com}\\ \\  E.\ Ventura \\
\small{Dept.\ Mat.\ Apl.\ III,}\\  {\small Universitat  Polit$\grave{\text {e} }$cnica de Catalunya,}
\\ {\small  Manresa, Barcelona, Catalunya}
\\ \small{e-mail: enric.ventura@upc.edu} }

\begin{document}

\maketitle

\begin{abstract}
Let $H$ be a torsion-free $\delta$-hyperbolic group with respect to a finite generating set $S$. Let $a_1,\ldots ,a_n$
and $a_{1*},\ldots ,a_{n*}$ be elements of $H$ such that $a_{i*}$ is conjugate to $a_i$ for each
$i=1,\dots ,n$. Then, there is a uniform conjugator if and only if $W(a_{1*},\ldots ,a_{n*})$ is conjugate to
$W(a_1,\ldots ,a_n)$ for every word $W$ in $n$ variables and length up to a computable constant depending only on
$\delta$, $\sharp{S}$ and $\sum_{i=1}^n |a_i|$.

As a corollary, we deduce that there exists a computable constant $\mathcal{C}=\mathcal{C}(\delta, \sharp S)$ such
that, for any endomorphism $\varphi$ of $H$, if $\varphi(h)$ is conjugate to $h$ for every element $h\in H$ of length
up to $\mathcal {C}$, then $\varphi$ is an inner automorphism.

Another corollary is the following: if $H$ is a torsion-free conjugacy separable hyperbolic group, then $\text{\rm
Out}(H)$ is residually finite.

When particularizing the main result to the case of free groups, we obtain a solution for a mixed version of the
classical Whitehead's algorithm. 

We show also that the Whitehead
problem and the mixed Whitehead problem for torsion-free hyperbolic groups are equivalent. 
\end{abstract}

\bigskip

\section{Introduction}

Let $G$ be a group and $A$ be a subset of $G$. An endomorphism $\varphi$ of $G$ is called \emph{pointwise inner on} $A$
if the element $\varphi(g)$ is conjugate to $g$, for every $g\in A$. We call $\varphi$ \emph{pointwise inner} if it is
pointwise inner on $G$. The group of all pointwise inner automorphisms of $G$ is denoted by $\text{\rm Aut}_{\text{\rm
pi}}(G)$. Clearly, $\text{\rm Inn}(G)\unlhd \text{\rm Aut}_{\text{\rm pi}}(G)\unlhd \text{\rm Aut}(G)$.

There are groups admitting pointwise inner automorphisms which are not inner. For example, some finite groups
(see~\cite{R}), some torsion-free nilpotent groups (see~\cite{Segal}), some nilpotent Lie groups
(see~\cite{GW}), and direct products of such groups with arbitrary groups. The fact that some nilpotent Lie groups
admit such automorphisms was used in~\cite{GW} to construct isospectral but not isometric Riemannian manifolds.

On the other hand, for free nilpotent groups (see~\cite{E}), for free groups (see~\cite{Gr,Lu}), for non-trivial free products (see~\cite{N}), and for fundamental
groups of closed surfaces of negative Euler characteristic (see~\cite{BKZ}), all pointwise inner automorphisms are
indeed inner. In the last paper, this property was used to show that surface groups satisfy a weak Magnus property.

One of the results in the present paper states that torsion-free hyperbolic groups also fall into this last class of
groups. In fact, we prove a stronger computational version of this fact: endomorphisms of torsion-free hyperbolic
groups which are pointwise inner on a ball of a uniformly bounded (and computable) radius, are indeed inner automorphisms.

\begin{thm}\label{cor}
Let $H$ be a torsion-free $\delta$-hyperbolic group with respect to a finite generating set $S$. Then, there exists a
computable constant $\mathcal{C}$ (depending only on $\delta$ and the cardinal $\sharp{S}$) such that, for every
endomorphism $\varphi$ of $H$, if $\varphi(g)$ is conjugate to $g$ for every element $g$ in the ball of radius
$\mathcal{C}$, then $\varphi$ is an inner automorphism.
\end{thm}

%\begin{cor}\label{algo}
%Given a finite presentation of a torsion-free $\delta$-hyperbolic group $H$, and an endomorphism %$\varphi$ of $H$ (by
%images of the generators), it is algorithmically decidable whether $\varphi$ is or is not an inner %automorphism.
%\end{cor}

%Clearly, a natural way of proving these results would be to find a uniform conjugator from a %previously given (big
%enough) set of pairs of conjugate elements. In fact, Theorem~\ref{cor} and Corollary~\ref{algo} %follow immediately from
%the main result in this paper:

An immediate consequence of Theorem~\ref{cor} is that one can algorithmically decide whether a given endomorphism of a
torsion-free hyperbolic group (given by a finite presentation, and images of generators) is or is not an inner
automorphism. This can also be easily deduced from the well-know fact that hyperbolic groups and their direct products
are bi-automatic; an alternative proof can also be found in~\cite[Theorem A]{BrHow}. However we stress, that the
purpose of the present paper is not the conjugacy problem for subsets of elements in hyperbolic groups.

\medskip

Theorem~\ref{cor} follows immediately from the main result of this paper:

\begin{thm}\label{maingenbeg}
Let $H$ be a torsion-free $\delta$-hyperbolic group with respect to a finite generating set $S$. Let $a_1,\ldots ,a_n$
and $a_{1*},\ldots ,a_{n*}$ be elements of $H$ such that $a_{i*}$ is conjugate to $a_i$ for every $i=1,\ldots ,n$.
Then, there is a uniform conjugator for them if and only if $W(a_{1*},\ldots ,a_{n*})$ is conjugate to $W(a_1,\ldots
,a_n)$ for every word $W$ in $n$ variables and length up to a computable constant depending only on $\delta$,
$\sharp{S}$ and $\sum_{i=1}^n |a_i|$.
\end{thm}

Note that Theorem~\ref{cor} was formulated in~\cite[Theorem 2]{BMV}. Independently, A.~Minasyan and D.~Osin~\cite{MO}
proved a variant of Theorem~\ref{maingenbeg}, for relatively hyperbolic groups but without the statement on
computability for the involved constant. Note also that our Theorem~\ref{cor} and~\cite[Theorem 1.1]{MO} both imply
that if $H$ is a torsion-free hyperbolic group, then the groups ${\text{\rm Inn}}(H)$ and ${\text{\rm Aut}}_{\text{\rm
pi}}(H)$ coincide.

V.~Metaftsis and M.~Sykiotis \cite{MS1, MS2} proved that, for any (relatively) hyperbolic group $H$, the group
${\text{\rm Inn}}(H)$ has finite index in ${\text{\rm Aut}}_{\text{\rm pi}}(H)$. Their proof is not constructive, it
uses ultrafilters and ideas of F.~Paulin on limits of group actions.

Furthermore, E.K.~Grossman proved in~\cite{Gr} that if $G$ is a finitely generated conjugacy
separable group, then the group ${\text{\rm Aut}}(G)/{\text{\rm Aut}}_{\text{\rm pi}}(G)$ is residually finite. From
this, one can immediately deduce the following corollary.

\begin{cor}\label{corGr} If $H$ is a torsion-free conjugacy separable hyperbolic group, then ${\text{\rm Out}}(H)$ is
residually finite.
\end{cor}

As a further application, we consider the case of a finitely generated free group $F$. Whitehead, back in 1936
(see~\cite{W} or~\cite{LS}), gave an algorithm to decide, given two tuples of elements of $F$, $a_1,\ldots ,a_n$ and
$b_1,\ldots ,b_n$, whether there is an automorphism of $F$ sending $a_i$ to a conjugate of $b_i$, for $i=1,\ldots ,n$
(with possibly different conjugators). Later, in 1974 (see~\cite{Mc} or~\cite{LS}), J.~McCool solved the same problem with
exact words: given two tuples of elements of $F$, $a_1,\ldots ,a_n$ and $b_1,\ldots ,b_n$, one can algorithmically
decide whether there is an automorphism of $F$ sending $a_i$ to $b_i$ for $i=1,\ldots ,n$. As a corollary of the main
result in the present paper, we obtain a mixed version of Whitehead/McCool's algorithm (see Theorem~\ref{Whitehead} for details). 

Finally we show, that if $H$ is a torsion-free hyperbolic group and the Whitehead problem in $H$ is solvable, then
the mixed Whitehead problem in $H$ is also solvable (see Theorem~\ref{WhiteheadHyp}). 

The structure of the paper is as follows. In Section~2 we recall some definitions and basic facts on hyperbolic metric
spaces and hyperbolic groups. Also, we prove there several statements (specially about norms and axes of elements, and
about controlling cancelations in some products of elements) which will be used later. The main theorem will be proved
in Sections~3 to~5, first in a very special case (Section~3), then in the case $n=2$ (Section~4), and finally in the
general case (Section~5). These three sections are sequential and the arguments in each one are helpful for the next
one. Finally, and particularizing the results to the case of free groups, in Section~6 we deduce a mixed version of
Whitehead's algorithm.

\section{Hyperbolic preliminaries}

\subsection{Hyperbolic spaces}

Let $(\mathcal{X}, d)$ be a metric space.

If $A,B$ are points or subsets of $\mathcal{X}$, the distance between them will be denoted by $d(A,B)$, or simply by
$|AB|$ if there is no risk of confusion.

A \emph{path} in $\mathcal{X}$ is a map $p\colon I\rightarrow \mathcal{X}$, where $I$ is an interval of the real line
(bounded or unbounded) or else the intersection of $\mathbb{Z}$ with such an interval. In the last case the path is
called \emph{discrete}. If $I=[a,b]$ then $p(a)$ and $p(b)$ are called the \emph{endpoints} of $p$. In that case we say
that the path $p$ is \emph{bounded} and \emph{goes from $p(a)$ to $p(b)$}; otherwise, we use the terms \emph{infinite
path} and \emph{bi-infinite path} with the obvious meaning. Sometimes we will identify a path with its image in
$\mathcal{X}$.

We say that a path $p$ is \emph{geodesic} if $d(p(r),p(s))=|r-s|$ for every $r,s\in I$. The space $(\mathcal{X}, d)$ is
said to be a \emph{geodesic} metric space if for every two points $A,B\in \mathcal{X}$ there is a geodesic from $A$ to
$B$ (not necessarily unique). Such a geodesic is usually denoted $[AB]$.

By a \emph{geodesic $n$-gon} $A_1A_2\cdots A_n$, where $n\geqslant 3$, we mean a cyclically ordered list of points
$A_1,\ldots ,A_n \in \mathcal{X}$ together with chosen geodesics $[A_1A_2], [A_2A_3],\ldots ,[A_{n-1}A_n], [A_nA_1]$;
each of these geodesics is called a \emph{side} of the $n$-gon, and each $A_i$ a \emph{vertex}. A geodesic 3-gon is
usually called a \emph{geodesic triangle}, and a geodesic 4-gon a \emph{geodesic rectangle}.

\begin{defn} {\rm
Let $(\mathcal{X},d)$ be a geodesic metric space and $\delta$ be a nonnegative real number.

A geodesic triangle $A_1A_2A_3$ in $\mathcal{X}$ is called \emph{$\delta$-thin} if for any vertex $A_i$ and any two
points $X\in [A_i,A_j]$, $Y\in [A_i,A_k]$ with
 $$
|A_iX|=|A_iY|\leqslant \frac{1}{2}(|A_iA_j|+|A_iA_k|-|A_jA_k|),
 $$
we have $|XY|\leqslant \delta$. The space $\mathcal{X}$ is called \emph{$\delta$-hyperbolic} if every geodesic
triangle in $\mathcal{X}$ is $\delta$-thin.
}\end{defn}

Directly from this definition it follows that each side of a $\delta$-thin triangle is contained in the
$\delta$-neighborhood of the union of the other two. By induction, one can easily extend this observation to $n$-gons.

\begin{prop}\label{gon}
If $A_1A_2\cdots A_n$ is a geodesic $n$-gon in a $\delta$-hyperbolic geodesic space, then each side is contained in the
$(n-2)\delta$-neighborhood of the union of all the others. \qed
\end{prop}

The following result is straightforward and will be used later (it is known as the \emph{rectangle inequality}).

\begin{prop}(see Remark 1.21 in~\cite[Chapter III.H]{BrH})\label{rec-in}
Any 4-gon $ABCD$ in a $\delta$-hyperbolic geodesic space $(\mathcal{X},d)$ satisfies the following inequality:
 $$
|AC|+|BD|\leqslant \max\{ |BC|+|AD|, |AB|+|CD|\}+2\delta. \,\,\,\Box
 $$
\end{prop}

Along the paper, we will need to use some approximations to the concept of geodesic. Here is a technical result and two
standard notions.

\begin{lem}\label{technical}
Let $A_1,A_2,\ldots , A_n$ be $n\geqslant 3$ points in a $\delta$-hyperbolic geodesic space satisfying the following
conditions:
\begin{itemize}
\item[(i)] $|A_{i-1}A_{i+1}|\geqslant |A_{i-1}A_i|+|A_iA_{i+1}|-2\delta$, for every $2\leqslant i\leqslant n-1$,
\item[(ii)] $|A_{i-1}A_i|>(2n-3)\delta$, for every $3\leqslant i\leqslant n-1$.
\end{itemize}
Then,
 \begin{equation}\label{eq16}
|A_1A_n|\geqslant \sum_{i=1}^{n-1}|A_iA_{i+1}|-(4n-10)\delta.
 \end{equation}
\end{lem}

\demo The proof goes by induction on $n$. Note that for $n=3$ the result is obvious.

Assume the result valid for $n$ points and let us prove it for $n+1$. Let $A_1,A_2,\ldots ,A_n,A_{n+1}$ be $n+1$ points
satisfying condition (i) for $2\leqslant i\leqslant n$, and condition (ii) for $3\leqslant i\leqslant n$. Clearly then
$A_1,A_2,\ldots ,A_n$ satisfy the corresponding conditions and, by the inductive hypothesis, we have
equation~(\ref{eq16}), so
 $$
|A_1A_n|\geqslant \sum_{i=1}^{n-1}|A_iA_{i+1}|-(4n-10)\delta \geqslant |A_1A_{n-1}|+|A_{n-1}A_n|-(4n-10)\delta.
 $$
From condition (i) with $i=n$ we have
 \begin{equation}\label{eq17}
|A_{n-1}A_{n+1}|\geqslant |A_{n-1}A_n|+|A_nA_{n+1}|-2\delta.
 \end{equation}
Adding these two last inequalities and applying condition (ii) for $i=n$, we get
 $$
|A_1A_n|+|A_{n-1}A_{n+1}|\geqslant |A_1A_{n-1}|+|A_nA_{n+1}|+2|A_{n-1}A_n|-(4n-8)\delta > |A_1A_{n-1}|+|A_nA_{n+1}|
+2\delta.
 $$
Therefore, the maximum in the rectangle inequality applied to $A_1A_{n-1}A_nA_{n+1}$ (see Proposition~\ref{rec-in}),
 $$
|A_1A_n|+|A_{n-1}A_{n+1}|\leqslant \max \{|A_1A_{n-1}|+|A_nA_{n+1}|,\, |A_1A_{n+1}|+|A_{n-1}A_n|\}+2\delta,
 $$
is achieved in the second entry. Hence,
 \begin{equation}\label{eq18}
|A_1A_n|+|A_{n-1}A_{n+1}|\leqslant |A_1A_{n+1}|+|A_{n-1}A_n|+2\delta.
 \end{equation}
On the other hand, from the induction hypothesis~(\ref{eq16}) and inequality~(\ref{eq17}), we have
 $$
\begin{array}{ll}|A_1A_n|+|A_{n-1}A_{n+1}| & \geqslant \Bigl(\sum_{i=1}^{n-1}|A_iA_{i+1}|-(4n-10)\delta\Bigr) +
|A_{n-1}A_n|+ |A_nA_{n+1}|-2\delta \\ \\ & =\sum_{i=1}^{n}|A_iA_{i+1}|+|A_{n-1}A_n|-(4n-8)\delta.
\end{array}
 $$
From this and inequality~(\ref{eq18}) we complete the proof:
 $$
|A_1A_{n+1}|\geqslant \sum_{i=1}^{n}|A_iA_{i+1}|-(4n-6)\delta = \sum_{i=1}^{n}|A_iA_{i+1}|-(4(n+1)-10)\delta. \quad
\Box
 $$

\begin{defn} \label{def} {\rm
Let $(\mathcal{X},d)$ be a metric space and $p\colon I\rightarrow \mathcal{X}$ be a path. Let $k>0$, $\lambda\geqslant
1$ and $\epsilon\geqslant 0$ be real numbers. The path $p$ is said to be \emph{$k$-local geodesic} if
$d(p(r),p(s))=|r-s|$ for all $r,s\in I$ with $|r-s|\leqslant k$. And it is said to be
\emph{$(\lambda,\epsilon)$-quasi-geodesic} if, for all $r,s\in I$, we have
 $$
\frac{1}{\lambda}|r-s|-\epsilon\leqslant d(p(r),p(s))\leqslant \lambda |r-s|+\epsilon.
 $$
}\end{defn}

\begin{prop}\label{1.13} (see Theorem 1.13~(3) in~\cite[Chapter~III.H]{BrH}).
Let $\mathcal{X}$ be a $\delta$-hyperbolic geodesic space and let $p\colon [a,b]\to \mathcal{X}$ be a $k$-local
geodesic with $k>8\delta$. Then, $p$ is a $(\lambda, \epsilon)$-quasi-geodesic, where
$\lambda=\frac{k+4\delta}{k-4\delta}$ and $\epsilon=2\delta$.
\end{prop}

The following proposition (without the statement on computability for $R$) is Theorem 1.7 in~\cite[Chapter III.H]{BrH}.
The computability of $R$ can be easily extracted from the proof there.

\begin{prop}\label{prop} (see Theorem 1.7 in~\cite[Chapter III.H]{BrH})
If $\mathcal{X}$ is a $\delta$-hyperbolic geodesic space, $p$ is a bounded $(\lambda,\epsilon)$-quasi-geodesic in
$\mathcal{X}$ and $c$ is a geodesic segment joining the endpoints of $p$, then {\rm im}\,$c$ and {\rm im}\,$p$ are
contained in the $R$-neighborhood of each other, where $R=R(\delta,\lambda,\epsilon)$ is a computable function.
\end{prop}

\subsection{Hyperbolic groups}

Let $H$ be a group given, together with a finite generating set $S$.

The length of an element $g\in H$ (with respect to $S$), denoted $|g|$, is defined as the length of the shortest word
in $S^{\pm 1}$ which equals $g$ in $H$. This naturally turns $H$ into a metric space; $|\, {\cdot}\, |$ is usually called the
\emph{word} metric.

Let $\Gamma(H,S)$ be the geometric realization of the right Cayley graph of $H$ with respect to $S$. We will consider
$\Gamma(H,S)$ as a metric space with the metric, induced by the word metric on $H$: $d(g_1,g_2)=|g_1^{-1}g_2|$. In
particular, edges are isometric to the real interval $[0,1]$. We highlight the fact that there is a notational
incoherence in using $|AB|$ to denote the distance between the points $A$ and $B$ in the Cayley graph $\Gamma(H,S)$,
while $|a^{-1}b|$ is the distance between the elements $a$ and $b$ of $H$; however, there will be no confusion because
we adopt the convention of using capital letters when thinking elements of $H$ as vertices of the Cayley graph.

The {\it ball of radius} $r$ around 1 in $\Gamma(H,S)$ is denoted $\mathcal{B}(r)$. The cardinality of any subset
$M\subseteq H$ is denoted $\sharp{M}$. For brevity, the cardinality of the set $\mathcal{B}(r)\cap H$ is denoted by
$\sharp \mathcal{B}(r)$. Clearly, an upper bound for $\sharp \mathcal{B}(r)$ is the number of elements in the similar
ball for the free group with basis $S$, so $\sharp \mathcal{B}(r)\leqslant 2(2\sharp S-1)^r$.

The group $H$ is called {\it $\delta$-hyperbolic with respect to $S$} if the corresponding metric space $\Gamma(H,S)$
is $\delta$-hyperbolic. It is well-known that if a group is hyperbolic with respect to some finite generating set, then
it is also hyperbolic with respect to any other finite generating set (with a possibly different $\delta$). This allows
to define hyperbolic groups: $H$ is said to be {\it hyperbolic} if for some finite generating set $S$, and some real
number $\delta\geqslant 0$, $H$ is $\delta$-hyperbolic with respect to $S$. It is also well-known that a finitely
generated group is free if and only if it is 0-hyperbolic with respect to some finite generating set $S$.

Let us begin with some well-known results about hyperbolic groups that will be needed later. The first one reproduces
Proposition 3.20 of~\cite[Chapter III.H]{BrH} plus the computability of the involved constant, which can be easily
extracted from the proof there. The second one solves the conjugacy problem within this family of groups. The following
one is about root-free elements in the torsion-free case ($g\in H$ is called {\it root-free} if it generates its own
centralizer, i.e. $C_H(g)=\langle g\rangle$). And the next one is also extracted from~\cite{BrH}.

\begin{prop}\label{free} (see Proposition~3.20 in~\cite[Chapter~III.H]{BrH}).
Let $H$ be a $\delta$-hyperbolic group with respect to a finite generating set $S$. For every finite set of elements
$h_1,\ldots ,h_r \in H$ there exists an integer $n>0$ such that $\langle h_1^n,\ldots ,h_r^n \rangle$ is free (of rank
$r$ or less). Furthermore, the integer $n$ is a computable function of $\delta$, $\sharp S$ and $\sum_i^r |h_i|$. \qed
\end{prop}

\begin{thm}\label{conj} (see Theorem 1.12 in~\cite[Chapter~III.$\Gamma$]{BrH}).
Let $H$ be a $\delta$-hyperbolic group with respect to a finite generating set $S$. If $u,v\in H$ are conjugate, then
the length of the shortest conjugator is bounded from above by a computable function of $\max\{|u|,|v|\}$, $\delta$ and
$\sharp S$. \qed
\end{thm}

\begin{lem}\label{minasyan} (see Lemma~4.3 in~\cite{Min})
Let $H$ be a torsion-free hyperbolic group, and let $a,b$ two elements, such that $b\notin C_H(a)$. Then there is a
computable integer $k_0=k_0(|a|,|b|)>0$, such that for every $k>k_0$ the element $ab^k$ is root-free. \qed
\end{lem}

\begin{prop}\label{quasigeod2} (see Corollary~3.10~(1) in~\cite[Chapter~III.$\Gamma$]{BrH}).
Let $H$ be a $\delta$-hyperbolic group with respect to a finite generating set $S$, and let $g\in H$ be an element of
infinite order. Then the map $\mathbb{Z}\rightarrow H$ given by $n\mapsto g^n$ is a quasi-geodesic. \qed
\end{prop}

The following lemma is well known and can be deduced straightforward from Proposition~\ref{quasigeod2}.

\begin{lem}\label{conjpowers}
Let $H$ be a $\delta$-hyperbolic group with respect to a finite generating set $S$, and let $g\in H$ be an element of
infinite order. If $g^p$ and $g^q$ are conjugate then $p=\pm q$. \qed
\end{lem}

Now, we provide an alternative proof for Proposition~\ref{quasigeod2}, in order to gain computability of the involved
constants.

\begin{lem}\label{computable}
The constants $\lambda$ and $\epsilon$ in Proposition~\ref{quasigeod2} are computable functions depending only on
$\delta$, $\sharp{S}$ and $|g|$.
\end{lem}

\demo First we make the following two easy observations:
\begin{itemize}
\item[(1)] Let $k\geqslant 1$ be a natural number and suppose that the map $\mathbb{Z}\rightarrow H$ given by $n\mapsto
g^{kn}$ is $(\lambda', \epsilon')$-quasi-geodesic. Then the map $\mathbb{Z}\rightarrow H$ given by $n\mapsto g^n$ is
$(\lambda,\epsilon)$-quasi-geodesic with $\lambda=k\lambda'$ and $\epsilon=\epsilon'+(k-1)|g|$. Thus, at any moment we
can replace $g$ by an appropriate power $g^k$.

\item[(2)] Let $g_0$ be a conjugate of $g$ in $H$, say $g=h^{-1}g_0h$ for some $h\in H$, and suppose that the map
$\mathbb{Z}\rightarrow H$, $n\mapsto g_0^n$, is $(\lambda', \epsilon')$-quasi-geodesic. Then, the map
$\mathbb{Z}\rightarrow H$, $n\mapsto g^n$, is $(\lambda, \epsilon)$-quasi-geodesic, where $\lambda=\lambda'$ and
$\epsilon=\epsilon'+2|h|$. Thus, at any moment we can replace $g$ by any conjugate $h^{-1}gh$.
\end{itemize}

Now, let us prove the result. Take an element $g\in H$ of infinite order. By Lemma~\ref{conjpowers},
there must exists an exponent $1\leqslant r\leqslant 1+\sharp \mathcal{B}(8\delta)$ such that the shortest conjugate of
$g^r$, say $g_0$, has length $|g_0|=k>8\delta$ (note that both $r$ and the corresponding conjugate are effectively
computable by Lemma~\ref{conj}). Replacing $g$ by $g_0$ and applying the previous two paragraphs, we may assume that
$|g|=k>8\delta$ and no conjugate of $g$ is shorter than $g$ itself.

Take a geodesic expression for $g$, say $g=s_1\cdots s_k$ with $s_i\in S^{\pm 1}$, and consider the bi-infinite path
$p_{g}\colon \mathbb{Z}\to H$ defined by the following rule: if $n\geqslant0$ and $n=tk+r$, where $0\leqslant r<k$,
then $p_g(n)=g^ts_1\cdots s_r$ and $p_g(-n)=g^{-t}s_k^{-1}\dots s^{-1}_{k-r+1}$; this corresponds to the bi-infinite
word $g^{\infty}=\cdots s_1\cdots s_ks_1\cdots s_k\cdots$. Clearly, any segment of length $k$ is of the form $s_i\cdots
s_ks_1\cdots s_{i-1}$, i.e. a conjugate of $g$ and hence geodesic. So, $p_g$ is a $k$-locall geodesic and thus a
$(8\delta+1)$-local geodesic. Finally, by Proposition~\ref{1.13}, $p_g$ is a $(3,2\delta)$-quasi-geodesic. Hence the
map $n\mapsto g^n$ is a $(3k,2\delta)$-quasi-geodesic. \qed

\medskip

Combining Proposition~\ref{prop} with Proposition~\ref{quasigeod2} and Lemma~\ref{computable}, we obtain the following
three corollaries.

\begin{cor} \label{useful}
Let $H$ be a $\delta$-hyperbolic group with respect to a finite generating set $S$, and let $g\in H$ be of infinite
order. Then for any integers $i<j$, the set $\{g^i, g^{i+1}, \ldots , g^j\}$ and any geodesic segment $[g^i,g^j]$
lie in the $\mu$-neighborhood of each other, where $\mu=\mu(\delta, \sharp S, |g|)$ is a computable function.
\end{cor}

\begin{cor}\label{powers}
Let $H$ be a $\delta$-hyperbolic group with respect to a finite generating set $S$, and let $g\in H$ be of infinite
order. For any natural numbers $s,t$ we have
 $$
|g^{s+t}|\geqslant |g^s|+|g^t|-2\mu,
 $$
where $\mu=\mu(\delta, \sharp S, |g|)$ is the constant from Corollary~\ref{useful}.
\end{cor}

\demo Consider the points $A=1$, $B=g^s$ and $C=g^{s+t}$ and choose geodesics $[AB]$, $[BC]$ and $[AC]$. By
Corollary~\ref{useful}, there exists $D\in [AC]$ such that $|BD|\leqslant\mu$. Then,
 $$
|AC|=|AD|+|DC|\geqslant (|AB|-|BD|)+(|CB|-|BD|)\geqslant |AB|+|BC|-2\mu. \quad \Box
 $$

The last corollary in this subsection is about torsion-free hyperbolic groups. It uses the following well known result.

\begin{prop}\label{nova}
Let $H$ be a torsion-free $\delta$-hyperbolic group. Then, centralizers of nontrivial elements are infinite cyclic. In
particular, extraction of roots is unique in $H$ (i.e. $g_1^r=g_2^r$ implies $g_1=g_2$). Furthermore, if for $1\neq
g\in H$, $g^p$ and $g^q$ are conjugate then $p=q$.
\end{prop}

\demo Cyclicity of centralizers is proven in~\cite[pages 462--463]{BrH}.

Suppose $g_1^r=g_2^r$. Then both $g_1$ and $g_2$ belong to the infinite cyclic group $C_H(g_1^r)$ and so, $g_1=g_2$.

Finally, suppose that $g^p=h^{-1}g^qh$; by Lemma~\ref{conjpowers}, $p=\epsilon q$ where $\epsilon =\pm 1$. Extracting
roots, $h^{-1}gh=g^{\epsilon}$. Thus, $h^2$ commutes with $g$ so both are powers of a common element, say $z\in H$. But
$h$ also commutes with $z$ so they are both powers of a common $y$, and so is $g$ too. Hence, $h^{-1}gh=g$ and
$\epsilon=1$. Thus, $p=q$. \qed

This proposition allows to use rational exponents in the notation, when working in torsion-free $\delta$-hyperbolic
groups (with $g^{1/s}$ meaning the unique element $x$ such that $x^s=g$, assuming it exists). For example, it is easy
to see that in such a group, every element commuting with $g^r\neq 1$ must be a rational power of $g$.

\begin{cor}\label{length}
Let $H$ be a torsion-free $\delta$-hyperbolic group with respect to a finite generating set~$S$. There exists a
computable function $f:\mathbb{N}^2\rightarrow \mathbb{N}$ such that, for any two elements $g,v\in H$ with $g$ of
infinite order, and for any nonnegative integers $p,q$ the following holds
 $$
|g^p vg^q|>|g^{p+q}|-f(|g|,|v|).
 $$
\end{cor}

\demo Let $\mu=\mu(|g|)$ be the computable constant given in Corollary~\ref{useful}: for any two integers $i<j$, the
set $\{ g^i,g^{i+1},\ldots ,g^j \}$ is contained in the $\mu$-neighborhood of any geodesic with endpoints $g^i$ and
$g^j$. Let $N=\sharp \mathcal{B}(2\delta+2\mu+|v|)$ and $M=2(N+1)(\mu+1)$.

Given $p,q\geqslant 0$, consider the points $A=1$, $B=g^p$, $C=g^p v$, and $D=g^p vg^q$, and choose geodesics $[AB]$,
$[AC]$, $[CD]$ and $[DA]$ (see Figure~1). Let $P$ be the point in $[CD]$ at distance $\ell=\frac{1}{2}(|AC|+|CD|-|AD|)$
from $C$.

\unitlength 1mm % = 2.85pt
\linethickness{0.4pt}
\ifx\plotpoint\undefined\newsavebox{\plotpoint}\fi % GNUPLOT compatibility
\begin{picture}(83.25,35.5)(-18,30)
\put(40,30){\line(1,0){40}} \put(51.5,55.5){\line(0,1){0}} \qbezier(40,30)(46.5,35)(52,55)
\qbezier(80,30)(73.5,35)(68,55) \qbezier(52,55)(60,53)(68,55) \qbezier(68,55)(58.5,39.5)(40,30)
\qbezier(40,30)(41.63,33.38)(39.75,36.25) \qbezier(80,30)(78.38,33.38)(80.25,36.25)
\qbezier(39.75,36.25)(42.25,38.88)(41.75,42) \qbezier(80.25,36.25)(77.75,38.88)(78.25,42)
\qbezier(41.75,42)(45.13,44.88)(44,48.25) \qbezier(78.25,42)(74.88,44.88)(76,48.25)
\qbezier(44,48.25)(46.63,49.25)(46.75,52.25) \qbezier(76,48.25)(73.38,49.25)(73.25,52.25)
\qbezier(46.75,52.25)(50,52)(52.25,54.75) \qbezier(73.25,52.25)(70,52)(67.75,54.75) \put(40,30){\circle*{1.5}}
\put(52,55){\circle*{1.5}} \put(68,55){\circle*{1.5}} \put(79.75,30.25){\circle*{1.5}} \put(73.25,39.5){\circle*{1.5}}
\put(71.5,44){\circle*{1.5}} \put(61.5,46.25){\circle*{1.5}} \put(50.25,48){\circle*{1.5}}
\qbezier(50.25,48.25)(55.25,46.75)(61.25,46.25) \qbezier(61.25,46.25)(66.13,43.88)(71.5,44)
\qbezier(50,48.25)(47,47.88)(44,48) \put(36.5,30.25){\makebox(0,0)[cc]{$A$}} \put(51.5,58.75){\makebox(0,0)[cc]{$B$}}
\put(68.5,58.75){\makebox(0,0)[cc]{$C$}} \put(83.25,30.25){\makebox(0,0)[cc]{$D$}}
\put(38.75,45.25){\makebox(0,0)[cc]{$g^p$}} \put(80.55,45.25){\makebox(0,0)[cc]{$g^q$}}
\put(59.75,56.75){\makebox(0,0)[cc]{$v$}} \put(46.25,38.25){\vector(1,2){.07}} \put(75.5,35.5){\vector(1,-2){.07}}
\put(58.5,54){\vector(1,0){3.25}} \put(71.,36.75){\makebox(0,0)[cc]{$P$}} \put(68.5,41.25){\makebox(0,0)[cc]{$X$}}
\put(43.25,50.5){\makebox(0,0)[cc]{$Y$}} \put(44.25,48){\circle*{1.5}} \put(77.75,50.5){\makebox(0,0)[cc]{$Z$}}
\put(75.75,48){\circle*{1.5}} \qbezier(71.5,44)(75,46)(75.75,48)
\end{picture}

\begin{center} Figure~1 \end{center}

If $\ell<M$ then
 $$
|g^p vg^q|=|AD|=|AC|+|CD|-2\ell\geqslant \bigl(|g^p|-|v|\bigr)+|g^q|-2\ell > |g^{p+q}|-|v|-2M.
 $$
Otherwise, if $\ell \geqslant M$ we will prove that $g$ and $v$ commute and so, $|g^p vg^q|=|g^{p+q}v|\geqslant
|g^{p+q}|-|v|$, concluding the proof.

So, assume $\ell \geqslant M$ and let us prove that $g$ and $v$ commute.

Let $X$ be an arbitrary point on $[CD]$ with $|CX|\leqslant \ell$. Then $X$ is at distance at most $\delta$ from the
side $[AC]$ of the geodesic triangle $ACD$. But this side is in the $(\delta+|v|)$-neighborhood of the side $[AB]$ of
the geodesic triangle $ABC$. And, by Corollary~\ref{useful}, this last one is in the $\mu$-neighborhood of the set
$\{1,g,\dots ,g^p\}$. Hence, there is a point of the form $Y=g^{p_0}$, $0\leqslant p_0\leqslant p$, such that
$|XY|\leqslant 2\delta+\mu+|v|$. Similarly, $X$ is in the $\mu$-neighborhood of $\{C,Cg,\ldots ,Cg^q\}$, i.e. there
exists a point of the form $Z=Cg^{q_0}=g^p vg^{q_0}$, $0\leqslant q_0\leqslant q$, such that $|XZ|\leqslant \mu$. Thus,
$|g^{p-p_0} vg^{q_0}|=|YZ|\leqslant |YX|+|XZ|\leqslant 2\delta+2\mu+|v|$.

Now, let $X_1,\ldots ,X_{N+1}$ be points on $[CD]$, such that $|CX_i|=2i(\mu+1)$ (the existence of all these points is
ensured by our assumption $\ell \geqslant M$). The previous paragraph gives us points $Y_i=g^{p_i}$ and $Z_i=g^p
vg^{q_i}$, with $0\leqslant p_i\leqslant p$ and $0\leqslant q_i\leqslant q$, such that $|X_iY_i|\leqslant
2\delta+\mu+|v|$ and $|X_iZ_i|\leqslant \mu$; thus, $|g^{p-p_i}vg^{q_i}|\leqslant 2\delta+2\mu+|v|$, for all
$i=1,\ldots ,N+1$. Furthermore, note that $q_i\neq q_j$ whenever $i\neq j$ (otherwise, $Z_i=Z_j$ and $|X_iX_j|\leqslant
|X_iZ_i|+|Z_jX_j|\leqslant 2\mu$, a contradiction).

This way we have obtained $N+1$ elements $g^{p-p_i}vg^{q_i}$ all of them in the ball $\mathcal{B}(2\delta+2\mu+|v|)$,
which has cardinal $N$. Thus, there must be at least one coincidence, $g^{p-p_i}vg^{q_i}=g^{p-p_j}vg^{q_j}$, for $i\neq
j$. Hence, $vg^{q_j-q_i}v^{-1}=g^{p_j-p_i}$. Since $q_i\neq q_j$, Proposition~\ref{nova} implies that $q_j-q_i=p_j-p_i$
and, extracting roots, $vgv^{-1}=g$. This means that $g$ commutes with $v$, completing the proof. \qed

\subsection{Controlling cancelation}

\begin{defn} {\rm
Let $H$ be a $\delta$-hyperbolic group with respect to a finite generating set $S$. For elements $u,v\in H$ and a real
number $c>0$ we write $uv=u\underset{{}^c}{\cdot}v$ if $\frac{1}{2}(|u|+|v|-|uv|)<c$. Also, we write
$uvw=u\underset{{}^c}{\cdot}v\underset{{}^c}{\cdot}w$ if $uv=u\underset{{}^c}{\cdot}v$ and $vw=v\underset{{}^c}{\cdot}w$. }\end{defn}

The definition of $u\underset{{}^c}{\cdot}v$ is equivalent to $|uv|>|u|+|v|-2c$. So, if $H$ is a free group,
$u\underset{{}^c}{\cdot}v$ means precisely that the maximal terminal segment of $u$ and the maximal initial segment of $v$
which can be canceled in the product $uv$ both have length smaller than $c$.

\begin{lem} \label{cc}
Let $H$ be a $\delta$-hyperbolic group with respect to a finite generating set $S$. If $c\in \mathbb{R}$ and $u,v,w\in
H$ are such that $uvw=u\underset{{}^c}{\cdot}v\underset{{}^c}{\cdot}w$ and $|v|>2c+\delta$, then
 $$
|u\underset{{}^c}{\cdot}v\underset{{}^c}{\cdot}w|> |u|+|v|+|w|-(4c+2\delta).
 $$
\end{lem}

\demo Connect the points $A=1$, $B=u$, $C=uv$ and $D=uvw$ by geodesic segments and consider the geodesic rectangle
$ABCD$. By assumption, $|BC|>2c+\delta$. From $u\underset{{}^c}{\cdot}v$ and $v\underset{{}^c}{\cdot}w$ we deduce
 $$
|AC|>|AB|+|BC|-2c>|AB|+\delta
 $$
and
 $$
|BD|>|BC|+|CD|-2c>|CD|+\delta,
 $$
respectively. From this and the rectangle inequality (Proposition~\ref{rec-in}), we deduce
 $$
(|AB|+|BC|-2c)+(|BC|+|CD|-2c)< |AC|+|BD|\leqslant |BC|+|AD|+2\delta,
 $$
which implies
 $$
|u|+|v|+|w|-(4c+2\delta )=|AB|+|BC|+|CD|-(4c+2\delta)<|AD|=|uvw|. \quad \Box
 $$

Next, we give some results about controlling cancelation that will be used later. Note that the important point in the
following lemma is the constant $c$ being independent from $k$.

\begin{lem}\label{circ}
Let $H$ be a $\delta$-hyperbolic group with respect to a finite generating set $S$, and let $w,b\in H$ with $b\neq 1$.
For every integer $k\geqslant 0$ and every $z\in H$, there exists $x\in H$ and $0\leqslant l\leqslant k$, such that
$z^{-1}wb^kz =x^{-1}\underset{{}^c}{\cdot} b^{k-l}wb^l \underset{{}^c}{\cdot} x$, where $c=3\delta+\mu(|b|)+|w|+1$ (and $\mu$ is
the computable function given in Corollary~\ref{useful}).
\end{lem}

\demo Fix $k\geqslant 0$ and $z\in H$, and let $0\leqslant l\leqslant k$ and $x\in H$ be such that
$z^{-1}wb^kz=x^{-1}b^{k-l}wb^lx$, with the shortest possible length for $x$; we will prove that these $l$ and $x$
satisfy the conclusion of the lemma. Suppose they do not, i.e. suppose that either
$x^{-1}b^{k-l}wb^l=x^{-1}\underset{{}^c}{\cdot}b^{k-l}wb^l$ or $b^{k-l}wb^lx=b^{k-l}wb^l \underset{{}^c}{\cdot}x$ is not true, and
let us find a contradiction. We consider only the case where the first of these expressions fails, i.e.
$|x^{-1}b^{k-l}wb^l|\leqslant |x^{-1}|+|b^{k-l}wb^l|-2c$; the second case can be treated analogously.

Consider the points $A=1$, $B=x^{-1}$, $C=x^{-1}b^{k-l}$, $D=x^{-1}b^{k-l}w$, $E=x^{-1}b^{k-l}wb^l$ and
$F=x^{-1}b^{k-l}wb^lx$, and connect them by geodesic segments, forming a 6-gon. In terms of the geodesic triangle
$ABE$, our assumption says $\frac{1}{2}(|AB|+|BE|-|AE|)\geqslant c$. By $\delta$-hyperbolicity of $H$, there exist
points $X_1\in [AB]$ and $X_2\in [BE]$ such that $|BX_1|=|BX_2|=c$ and $|X_1X_2|\leqslant \delta$. And, by
Proposition~\ref{gon} applied to the rectangle $BCDE$, there exists a point $X_3\in [BC]\cup [CD]\cup [DE]$ such that
$|X_2X_3|\leqslant 2\delta$.

\medskip

\emph{Case 1:} $X_3\in [BC]$ (see Figure~2). Since $C=Bb^{k-l}$, Corollary~\ref{useful} implies that there exists an
element $X_4=Bb^s$ for some $0\leqslant s\leqslant k-l$, such that $|X_3X_4|\leqslant \mu(|b|)$. Hence,
$|X_1X_4|\leqslant |X_1X_2|+|X_2X_3|+|X_3X_4|\leqslant 3\delta+\mu(|b|)<c$ and
$z^{-1}wb^kz=X_4b^{k-l-s}wb^{l+s}X_4^{-1}$.

\medskip

\emph{Case 2:} $X_3\in [CD]$. In this case, take $X_4=C$ and we have $|X_1X_4|\leqslant
|X_1X_2|+|X_2X_3|+|X_3X_4|\leqslant 3\delta+|w|<c$ as well. Similarly, $z^{-1}wb^kz =X_4wb^{k}X_4^{-1}$.

\medskip

\emph{Case 3:} $X_3\in [DE]$. Since $E=Db^l$, Corollary~\ref{useful} implies again that there exist an element
$X_4=Db^s$ for some $0\leqslant s\leqslant l$, such that $|X_3X_4|\leqslant \mu(|b|)$. Like in Case~1, we have
$|X_1X_4|<c$ and $z^{-1}wb^kz =X_4b^{k-s}wb^sX_4^{-1}$.

\medskip

In any case, we have found an element $X_4 \in H$ and a decomposition of $z^{-1}wb^kz$ of the form
$z^{-1}wb^kz=X_4b^{k-s}wb^sX_4^{-1}$, with $0\leqslant s\leqslant k$ and $|X_1X_4|<c$. Since $|X_1B|=c$, we have
 $$
|X_4|=|AX_4|\leqslant |AX_1|+|X_1X_4|<|AX_1|+|X_1B|=|AB|=|x|,
 $$
contradicting the minimality of $|x|$. \qed

\unitlength 1mm % = 2.85pt
\linethickness{0.4pt}
\ifx\plotpoint\undefined\newsavebox{\plotpoint}\fi % GNUPLOT compatibility
\begin{picture}(102.25,40.5)(-8,26)
\put(40,30){\line(1,0){60}} \qbezier(40,30)(45,32)(47,50) \qbezier(100,30)(95,32)(93,50)
\qbezier(47,50)(70,37.5)(93,50) \qbezier(47,50)(59.5,49)(64,56) \qbezier(93,50)(80.5,49)(76,56)
\qbezier(93,50)(80.5,49)(76,56) \qbezier(64,56)(70,50)(76,56) \qbezier(47,50)(50.5,50.75)(51,53.5)
\qbezier(93,50)(89.5,50.75)(89,53.5) \qbezier(51.25,53.25)(53.5,53.13)(54.75,55.5)
\qbezier(88.75,53.25)(86.5,53.13)(85.25,55.5) \qbezier(54.75,55.25)(57.25,54.63)(58.75,56.5)
\qbezier(85.25,55.25)(82.75,54.63)(81.25,56.5) \qbezier(58.75,56.5)(61.63,54.25)(64,56)
\qbezier(81.25,56.5)(78.38,54.25)(76,56) \put(40,30){\circle*{1.25}} \put(100,30){\circle*{1.25}}
\put(47,50){\circle*{1.25}} \put(93,50){\circle*{1.25}} \put(64,56){\circle*{1.25}} \put(76,55.75){\circle*{1.25}}
\put(45.50,40.00){\circle*{1}} \put(59.5,51.75){\circle*{1}} \put(58,45.50){\circle*{1}} \put(58.75,56.25){\circle*{1}}
\put(37.00,30.25){\makebox(0,0)[cc]{$A$}} \put(43.75,51){\makebox(0,0)[cc]{$B$}}
\put(64.75,59.00){\makebox(0,0)[cc]{$C$}} \put(76.00,58.75){\makebox(0,0)[cc]{$D$}}
\put(96.00,51.00){\makebox(0,0)[cc]{$E$}} \put(103.,30.25){\makebox(0,0)[cc]{$F$}}
\put(58.5,42.5){\makebox(0,0)[cc]{$X_2$}} \put(62.25,49.75){\makebox(0,0)[cc]{$X_3$}}
\put(58.5,59.00){\makebox(0,0)[cc]{$X_4$}} \put(69.75,55.5){\makebox(0,0)[cc]{$w$}}
\put(42.,36.5){\makebox(0,0)[cc]{$x$}} \put(98,37){\makebox(0,0)[cc]{$x$}} \put(48.75,56){\makebox(0,0)[cc]{$b^{k-l}$}}
\put(85.25,57.5){\makebox(0,0)[cc]{$b^l$}} \put(42.60,32.5){\vector(-2,-3){.07}} \put(97.5,32.5){\vector(1,-2){.07}}
\put(54,54.25){\vector(3,2){.07}} \put(68.75,53){\vector(1,0){2.5}} \put(87.75,53.25){\vector(3,-2){.07}}
\multiput(58.75,56.75)(.032609,-.228261){23}{\line(0,-1){.228261}}
\multiput(59.5,51.5)(-.044398947,-.1907895){33}{\line(0,-1){.1907895}}
\multiput(57.75,45.25)(-.08053691,-.03355705){149}{\line(-1,0){.08053691}} \put(48,38){\makebox(0,0)[cc]{$X_1$}}
\end{picture}

\begin{center} Figure 2 \end{center}

The previous lemma in the particular case of $w=1$ says that, for every $b,z\in H$ and every $k\geqslant 0$, there
exists $x\in H$ such that $z^{-1}b^kz =x^{-1}\underset{{}^c}{\cdot} b^k \underset{{}^c}{\cdot} x$, (where $c$ is a computable
function depending only on $\delta$ and $|b|$). In the following result we present a technical improvement (which will
be crucial later) showing that, in fact, one can choose a uniform $x$ valid for every $k$.

\begin{lem}\label{easy1}
Let $H$ be a $\delta$-hyperbolic group with respect to a finite generating set $S$, and let $z,b\in H$. There exists an
element $x\in H$ such that for every integer~$k$ holds $z^{-1}b^kz=x^{-1}\underset{{}^c}{\cdot}b^k\underset{{}^c}{\cdot}x$, where
$c=\delta+\mu(|b|)$.
\end{lem}

\demo Let $x^{-1}$ be one of the shortest elements in the set $\mathcal{G}=\{z^{-1}b^n\,|\, n\in \Z\}$. Clearly
$z^{-1}b^kz=x^{-1}b^kx$ for every $k\in \mathbb{Z}$. We show that $z^{-1}b^kz=x^{-1}\underset{{}^c}{\cdot}
b^k\underset{{}^c}{\cdot}x$. Fix $k\in \mathbb{Z}$ and denote $A=1$, $B=x^{-1}$, and $C=x^{-1}b^k$. We choose geodesic
segments $[AB]$, $[BC]$ and $[AC]$ and consider the points $X\in [BA]$, $Y\in [BC]$ such that $|BX|=|BY|=
\frac{1}{2}(|BA|+|BC|-|AC|)$. By $\delta$-hyperbolicity we have $|XY|\leqslant \delta$. By Corollary~2.14, the point
$Y\in [BC]$ lies at distance at most $\mu(|b|)$ from a point $D\in \mathcal{G}$. By the choice of $x^{-1}$, we have
$|AB|\leqslant |AD|$ and so
 $$
|AX|+|XB|=|AB|\leqslant |AD|\leqslant |AX|+|XY|+|YD|\leqslant |AX|+\delta+\mu(|b|).
 $$
Hence $|XB|\leqslant c$, i.e. $\frac{1}{2}(|x^{-1}|+|b^k|-|x^{-1}b^k|)\leqslant c$ and hence, $x^{-1}b^k=x^{-1}
\underset{{}^c}{\cdot}b^k$. Inverting the last element, and changing $k$ by $-k$, we have $b^{k}x=b^k\underset{{}^c}{\cdot}x$.
Thus, $x^{-1}b^{k}x=x^{-1}\underset{{}^c}{\cdot}b^k\underset{{}^c}{\cdot}x$. \qed

\subsection{The norm and the axis of an element}

\begin{defn} {\rm
Let $H$ be a torsion-free $\delta$-hyperbolic group with respect to a finite generating set $S$, and let $g\in H$. The
\emph{norm of $g$}, denoted $||g||$, is defined as
 $$
\min \{d(x,gx) \mid x\in \Gamma(H,S) \}.
 $$
The \emph{axis of $g$}, denoted $\mathcal{A}_g$, is the set of points $x\in \Gamma(H,S)$ where this minimum is
achieved,
 $$
\mathcal{A}_g =\{ x\in \Gamma(H,S)\,\mid \, d(x,gx)=||g|| \}.
 $$
}\end{defn}

The following facts are easy to see:
\begin{itemize}
\item[(1)] $\mathcal{A}_g\cap H$ is nonempty, in particular
$$
||g||=\min \{|x^{-1}gx| \mid x\in H \}.
$$
Moreover, $\mathcal{A}_g$ lies in the 1-neighborhood of $\mathcal{A}_g\cap H$;
\item[(2)] $||g||$ is a nonnegative integer satisfying $0\leqslant ||g||\leqslant |g|$. Moreover, $||g||=0$ iff $g=1$;
\item[(3)] $\mathcal{A}_g$ is $C_H(g)$-invariant: for every $x\in \mathcal{A}_g$ and $h\in C_H(g)$ we have $hx\in
\mathcal{A}_g$;
\item[(4)] for any $x\in \mathcal{A}_g$, any geodesic segment $[x,gx]$ also lies in $\mathcal{A}_g$;
\item[(5)] for any $h\in H$ we have $||hgh^{-1}||=||g||$ and $\mathcal{A}_{hgh^{-1}}=h\mathcal{A}_g$;
\item[(6)] for any $g\in H$ and any $x\in \Gamma(H,S)$, we have $d(x,gx)\leqslant ||g||+2d(x,\mathcal{A}_g)$.
\end{itemize}

\begin{lem}\label{noteasy}
Let $H$ be a torsion-free $\delta$-hyperbolic group with respect to a finite generating set $S$. For any $1\neq g\in
H$, there exists a computable integer $r=r(|g|)\geqslant 1$ such that
 $$
\bigcup_{k=1}^{\infty}\mathcal{A}_{g^k} \subseteq \langle g \rangle \mathcal{B}(r).
 $$
\end{lem}

\demo By Property (1), $\bigcup_{k=1}^{\infty}\mathcal{A}_{g^k}$ lies in the 1-neighborhood of $\bigcup_{k=1}^{\infty}
\mathcal{A}_{g^k}\cap H$. The strategy now is to see that this last set lies at bounded (in terms of $|g|$) distance
from the centralizer $C_H(g)$; and then, we will see that $C_H(g)$ lies at bounded distance from $\langle g\rangle$.

Take an arbitrary $z\in \cup_{k=1}^{\infty} \mathcal{A}_{g^k}\cap H$. By Properties~(1)-(2), there is $k\geqslant 1$
such that $|z^{-1}g^k z|$ is minimal among the lengths of all conjugates of $g^k$ (in particular, $|z^{-1}g^k
z|\leqslant |g^k|$). By Corollary~\ref{easy1}, there exists $x\in H$ such that $z^{-1}g^k z=x^{-1}\underset{{}^c}{\cdot} g^k
\underset{{}^c}{\cdot} x$, where the constant $c=c(|g|)$ is computable and independent from $k$. Thus, we have
$|x^{-1}\underset{{}^c}{\cdot} g^k \underset{{}^c}{\cdot} x|\leqslant |g^k|$. Let us consider two cases.

\medskip

{\it Case 1: $|g^k|> 2c+\delta$}. By Lemma~\ref{cc}, $|x^{-1}\underset{{}^c}{\cdot} g^k \underset{{}^c}{\cdot} x|>
2|x|+|g^k|-(4c+2\delta)$. Therefore $|x|< 2c+\delta$. Moreover, $z\in C_H(g)x$.

{\it Case 2: $|g^k|\leqslant 2c+\delta$}. From $|z^{-1}g^kz|\leqslant |g^k|$ and Theorem~\ref{conj}, we conclude that
there exists $y\in H$ such that $z^{-1}g^kz=y^{-1}g^ky$ and the length of $y$ is bounded by a computable constant,
depending only on $|g|$ (i.e. on $2c+\delta$). Moreover, $z\in C_H(g)y$.

\medskip

In both cases $z$ lies at bounded (in terms of $|g|$) distance from $C_H(g)$.

It remains to prove that $C_H(g)$ is at bounded distance from $\langle g\rangle$. Let $z\in C_H(g)$. By
Lemma~\ref{conjpowers}, there exists a (computable) natural number $s\leqslant \sharp \mathcal{B}(4\delta)$, such that
$g^s$ is not conjugate into the ball $\mathcal{B}(4\delta)$. In this situation, the proof of Corollary~3.10
in~\cite[Chapter~III.$\Gamma$]{BrH} shows that the distance from $z$ to the set $\langle g^s\rangle$ is at most
$2|g^s|+4\delta$. Hence, the distance from $z$ to $\langle g\rangle$ is bounded by a computable constant depending only
on $\delta, \sharp S$ and $|g|$.\qed

\medskip

From this lemma, it is easy to deduce the following corollaries.

\begin{cor}\label{0}
Let $H$ be a torsion-free $\delta$-hyperbolic group with respect to a finite generating set $S$. For any $1\neq g\in H$
and any integer $k\neq 0$, there exists an element $x\in \mathcal{A}_{g^k}\cap H$ of length at most $r(|g|)$. \qed
\end{cor}

\begin{cor}\label{1}
Let $H$ be a torsion-free $\delta$-hyperbolic group with respect to a finite generating set $S$. For any $1\neq g\in H$
and any integer $k\neq 0$, we have $||g^k||\geqslant |g^k|-2r(|g|)$.
\end{cor}

\demo Take the element $x$ from Corollary~\ref{0}. Then $||g^k||=d(x,g^kx)=|x^{-1}g^kx|\geqslant |g^k|-2|x|\geqslant
|g^k|-2r(|g|)$. \qed

\begin{cor}\label{2}
Let $H$ be a torsion-free $\delta$-hyperbolic group with respect to a finite generating set $S$. For any $1\neq g\in H$
and any $C>0$, there exists a computable integer $k_0=k_0(|g|,C)$ such that for any $k>k_0$ we have $||g^k||>C$.
\end{cor}

\demo Using Corollary~\ref{1}, and Proposition~\ref{quasigeod2} complemented with Lemma~\ref{computable}, we deduce
$||g^k||\geqslant |g^k|-2r(|g|)\geqslant \frac{1}{\lambda}k-\epsilon-2r(|g|)$ for every $k>0$, where $\lambda$,
$\epsilon$ and $r$ are computable functions of $|g|$. Now, the result follows easily. \qed

\begin{cor}\label{transl}
Let $H$ be a torsion-free $\delta$-hyperbolic group with respect to a finite generating set $S$. There exist computable
functions $f_1:\mathbb{N}\rightarrow \mathbb{N}$ and $f_2:\mathbb{N}\rightarrow \mathbb{N}$ such that, for every $1\neq
g\in H$ and every natural numbers $s,t>0$, we have
 $$
||g^{s+t}||-f_1(|g|)\leqslant ||g^s||+||g^t||\leqslant ||g^{s+t}||+f_2(|g|).
 $$
\end{cor}

\demo Take $f_1(n)=4r(n)$ and the first inequality follows from Corollary~\ref{1}:
 $$
||g^{s+t}||\leqslant |g^{s+t}|\leqslant |g^s|+|g^t|\leqslant ||g^s||+||g^t||+4r(|g|).
 $$
And taking $f_2(n)=2r(n)+2\mu(n)$, the second inequality follows from Corollaries~\ref{1} and~\ref{powers}:
 $$
||g^s||+||g^t||\leqslant |g^s|+|g^t|\leqslant |g^{s+t}|+2\mu(|g|)\leqslant ||g^{s+t}||+2r(|g|)+2\mu(|g|). \quad \Box
 $$

\medskip

Next, we will state several lemmas about distances to axes.

\begin{lem}\label{6.1}
Let $H$ be a torsion-free $\delta$-hyperbolic group with respect to a finite generating set $S$. Let $1\neq g\in H$,
let $A$ be a point in $\Gamma(H,S)$, and let $B$ be a point in $\mathcal{A}_{g}$ at minimal distance from $A$. Then,
for every geodesic segment $[BC]\subset \mathcal{A}_{g}$, we have
 $$
|AC|\geqslant |AB|+|BC|-2\delta.
 $$
\end{lem}

\demo Consider a given geodesic segment $[BC]$ contained in $\mathcal{A}_g$, and choose geodesic segments $[AB]$ and
$[AC]$. Let $X\in [BA]$ and $Y\in [BC]$ be points such that $|BX|=|BY|=\frac{1}{2}(|BA|+|BC|-|AC|)$. Then
$|XY|\leqslant \delta$. Since the point $Y$ also lies on $\mathcal{A}_g$, we have that $|AB|\leqslant |AY|$. Therefore
$|XB|\leqslant |XY|\leqslant \delta$. Thus,
 $$
|AC|=|AB|+|BC|-2|BX|\geqslant |AB|+|BC|-2\delta. \quad \Box
 $$

\medskip

\begin{lem} \label{5.7}
Let $H$ be a torsion-free $\delta$-hyperbolic group with respect to a finite generating set $S$. Let $g\in H$, and let
$k$ be an integer number such that $||g^k||>5\delta$. Let $A$ be an element of $H$, and $n\geqslant 0$ be such
that $d(A,g^kA)=||g^k||+n$. Then, $A=g^tv$ for some $t\in \mathbb{Z}$ and $v\in H$ with $|v|\leqslant
\frac{n}{2}+3\delta+r(|g|)$, where $r$ is the function introduced in Lemma~\ref{noteasy}.
\end{lem}

\demo By the hypothesis, $g\neq 1$. Let $B$ be a point in $\mathcal{A}_{g^k}$ at minimal distance from $A$. Let
$C=g^kB$ and $D=g^kA$. Since $C\in \mathcal{A}_g$ is at minimal distance from $D$ (the same as $|AB|$), Lemma \ref{6.1}
tells us that
 $$
|AC|\geqslant |AB|+|BC|-2\delta
 $$
and
 $$
|DB|\geqslant |CD|+|BC|-2\delta.
 $$
Moreover, $|BC|=||g^k||>5\delta$. Therefore, by Lemma~\ref{technical} applied to points $A,B,C,D$, we deduce
 $$
\begin{array}{rcl} |AD| & \geqslant & |AB|+|BC|+|CD|-6\delta \\ & = & 2|AB|+||g^k||-6\delta.
\end{array}
 $$
Hence, $|AB|\leqslant \frac{n}{2}+3\delta$. By Lemma~\ref{noteasy}, $B$ lies at distance at most $r(|g|)$ from $\langle
g\rangle$. Hence, $A$ lies at distance at most $\frac{n}{2}+3\delta+r(|g|)$ from $\langle g\rangle$. This completes the
proof. \qed

\begin{lem} \label{6.3}
Let $H$ be a torsion-free $\delta$-hyperbolic group with respect to a finite generating set $S$, and let $g\in H$ with
$||g||>5\delta$. Then the middle point of any geodesic segment $[A,gA]$, where $A$ is a point of $\Gamma(H,S)$, lies in
the $(5\delta)$-neighborhood of the axis $\mathcal{A}_g$.
\end{lem}

\demo By the hypothesis, $g\neq 1$. Let $B$ be a point in $\mathcal{A}_g$ at minimal distance from $A$. Let $C=gB$ and
$D=gA$. Exactly like in the previous lemma, we obtain
 \begin{equation}\label{19}
2|AB|+|BC|\leqslant |AD|+6\delta.
 \end{equation}
Now, take geodesic segments $[AD]$ and $[BC]$, and let $M$ and $N$ be their middle points, respectively. Clearly, $N\in
\mathcal{A}_g$. In order to estimate the distance $|NM|$, we consider the geodesic rectangle $AMDN$. By the rectangle
inequality, we have
 $$
\begin{array}{rcl} |NM|+|AD| & \leqslant & \max\{|AM|+|DN|,|DM|+|AN|\}+2\delta \\ & = & \max\{ \frac{1}{2}|AD|+|DN|,
\frac{1}{2}|AD|+|AN|\}+2\delta. \end{array}
 $$
But $|AN|\leqslant |AB|+|BN|=|AB|+\frac{1}{2}|BC|$. Therefore from~(\ref{19}), we have $|AN|\leqslant
\frac{1}{2}|AD|+3\delta$. Analogously, $|DN|\leqslant \frac{1}{2}|AD|+3\delta$. From all this we deduce
 $$
|NM|+|AD|\leqslant \frac{1}{2}|AD|+\frac{1}{2}|AD|+3\delta+2\delta.
 $$
Thus, $|NM|\leqslant 5\delta$. \qed

\begin{prop}\label{axes}
Let $H$ be a torsion-free $\delta$-hyperbolic group with respect to a finite generating set $S$, and let $g,h\in H$
with $||g||>15\delta$, $||h||>15\delta$ and $||gh||>5\delta$. Then the distance between the axes $\mathcal{A}_g$ and
$\mathcal{A}_h$ is at most
 $$
\max \{15\delta,\, \frac{1}{2}(||gh||-||g||-||h||)+18\delta\}.
 $$
\end{prop}

\demo By the hypotheses, $g,h$ and $gh$ are all nontrivial. Let $d=d(\mathcal{A}_g, \mathcal{A}_h)$, and let $X\in
\mathcal{A}_h$ and $Y\in \mathcal{A}_g$ be such that $|XY|=d$. If $d\leqslant 15\delta$ we are done so, let us assume
$d>15\delta$.

Consider the points $A_1=X$, $A_2=Y$, $A_3=gY$, $A_4=gX$, $A_5=ghX$, $A_6=ghY$, $A_7=ghgY$, $A_8=ghgX$, and
$A_9=ghghX$. By Lemma~\ref{6.1} and doing the appropriate translation, we have $|A_{i-1}A_{i+1}|\geqslant
|A_{i-1}A_i|+|A_iA_{i+1}|-2\delta$ for every $i=2,\ldots ,8$. Moreover, $|A_{i-1}A_i|$ equals either $d$, or $||g||$,
or $||h||$ which are all bigger than $15\delta$. So, Lemma~\ref{technical} tells us that
 $$
\begin{array}{rcl} d(A_1, A_9)=d(X,(gh)^2X) & \geqslant & d(X,Y)+d(Y,gY)+d(gY,gX)+d(gX,ghX)+d(ghX,ghY) \\ & \\ & &
+d(ghY,ghgY)+d(ghgY,ghgX)+d(ghgX,ghghX)-26\delta \\ & \\ & = & 2(d+||g||+d+||h||)-26\delta. \end{array}
 $$
On the other hand,
 $$
d(A_1,A_5)=d(X,ghX)\leqslant d(X,Y)+d(Y,gY)+d(gY,gX)+d(gX,ghX)=d+||g||+d+||h||.
 $$
Let now $[A_1A_5]$ be a geodesic segment, and consider its translation $(gh)[A_1A_5]$, say $[A_5A_9]$. Let $M$ be the
middle point of $[A_1A_5]$ and $M'=ghM$ be the middle point of $[A_5A_9]$. Since $\frac{1}{2}d(A_1,A_5)=d(A_1,M)
=d(M,A_5)=d(M',A_9)$, using the previous inequalities we have
 $$
\begin{array}{rcl} d(M,M') & \geqslant & d(A_1,A_9)-d(A_1,M)-d(M',A_9) \\ & = & d(A_1,A_9) -d(A_1,A_5) \\ &
\geqslant & 2d+||g||+||h||-26\delta. \end{array}
 $$
Finally, by Lemma~\ref{6.3}, $M$ lies at distance at most $5\delta$ from the axis $\mathcal{A}_{gh}$. Therefore,
$d(M,M')=d(M,ghM)\leqslant 10\delta+||gh||$. Hence $d\leqslant \frac{1}{2}(||gh||-||g||-||h||)+18\delta.$ $\Box$

\bigskip

\section{A special case of the main Theorem}

In this section, we prove a special case of Theorem~\ref{maingenbeg}, namely the case of two words ($n=2$) and with the
extra assumption that $\langle a_1, a_2\rangle$ is a cyclic subgroup of $H$. The proof contains ingredients which will
be used for the general case.

Let us start with the following lemma, which considers the situation where the product of conjugates of two powers of a
given element equals the product of these powers, and analyzes how the involved conjugators must look like.

\begin{lem}\label{5}
Let $H$ be a torsion-free $\delta$-hyperbolic group with respect to a finite generating set $S$. There exists a
computable function $\hbar:\mathbb{N}\rightarrow \mathbb{R}^{+}$ with the following property: for any three elements
$b,x,y\in H$ and any two positive integers $s,t$, which satisfy\,\, $||b^s||, ||b^t||>15\delta$, $||b^{s+t}||>5\delta$
and
 \begin{equation}\label{eq13}
(x\cdot b^s \cdot x^{-1})(y\cdot b^t \cdot y^{-1})=b^{s+t},
 \end{equation}
there exist integers $n_1,n_2,n_3,n_4$ and elements $v_x,v_y\in H$ of length at most $\hbar(|b|)$ such that
 $$
x=b^{n_1}v_x b^{n_2} \quad \quad \text{and} \quad \quad y=b^{n_3}v_y b^{n_4}.
 $$
\end{lem}

\demo Let $b,x,y$ and $s,t$ be as in the statement (in particular, $b\neq 1$). Consider the axes $\mathcal{A}_{xb^s
x^{-1}}=x\mathcal{A}_{b^s}$ and $\mathcal{A}_{yb^t y^{-1}}=y\mathcal{A}_{b^t}$. By Proposition~\ref{axes} applied to
the elements $xb^s x^{-1}$ and $yb^t y^{-1}$ (note that $||xb^s x^{-1}||=||b^s||>15\delta$, $||yb^t
y^{-1}||=||b^t||>15\delta$ and $||(xb^s x^{-1})(yb^t y^{-1})||=||b^{s+t}||>5\delta$ by hypothesis), the distance
between $x\mathcal{A}_{b^s}$ and $y\mathcal{A}_{b^t}$ is at most
 $$
\max \{15\delta,\, \frac{1}{2}(||b^{s+t}||-||b^s||-||b^t||)+18\delta\}.
 $$
By Corollary~\ref{transl}, this value does not exceed $\frac{1}{2}f_1(|b|)+18\delta$, an upper bound which is
independent from $s$ and $t$.

Now, take an element $Q\in y\mathcal{A}_{b^t}\cap H$ such that $d(Q, x\mathcal{A}_{b^s})\leqslant \frac{1}{2}f_1(|b|)
+18\delta+1$, and set $P=(yb^t y^{-1})^{-1}Q$. In particular, $P\in y\mathcal{A}_{b^t}\cap H$ and
$d(P,Q)=||yb^ty^{-1}||=||b^t||$. Then we have
 $$
\begin{array}{ll} d(P, b^{s+t}P) & =d(P, (xb^s x^{-1})(yb^t y^{-1})P)= d(P, (xb^s x^{-1})Q) \\ \\ & \leqslant d(P,Q)+d(Q,
(xb^s x^{-1})Q) \\ \\ & \leqslant d(P,Q)+2d(Q, \mathcal{A}_{xb^sx^{-1}})+||b^s||\\ \\ & \leqslant ||b^t||+||b^s||
+f_1(|b|)+36\delta +2\\ \\ & \leqslant ||b^{s+t}||+f_1(|b|)+f_2(|b|)+36\delta+2, \end{array}
 $$
where the last inequality uses Corollary~\ref{transl} again. Next, apply Lemma~\ref{5.7} to conclude that
$P=b^{n_3}v_1$ for some $n_3\in \mathbb{Z}$ and $v_1\in H$ with $|v_1 |\leqslant \frac{1}{2}f_1(|b|)
+\frac{1}{2}f_2(|b|)+r(|b|)+21\delta+1$. And since $P\in y\mathcal{A}_{b^t}\cap H$, we deduce from Lemma~\ref{noteasy}
that $y^{-1}P=b^{-n_4}v_2$, for some $n_4\in \mathbb{Z}$ and $v_2 \in H$ with $|v_2 |\leqslant r(|b|)$. Hence,
 $$
y=b^{n_3}v_y b^{n_4},
 $$
where $v_y =v_1v_2^{-1}$ has length bounded by
 $$
|v_y |=|v_1v_2^{-1}|\leqslant |v_1|+|v_2| \leqslant \frac{1}{2}f_1(|b|)+ \frac{1}{2}f_2(|b|)+2r(|b|)+21\delta+1.
 $$

Finally, inverting and replacing $b$ to $b^{-1}$ in equation~(\ref{eq13}), we obtain again the same equation with $x$
and $y$ interchanged. So, the same argument shows that
 $$
x=b^{n_1}v_x b^{n_2},
 $$
for some $n_1,n_2\in \mathbb{Z}$ and some $v_x\in H$ with the same upper bound for its length.

Hence, the function $\hbar(n)=\frac{1}{2}f_1(n) +\frac{1}{2}f_2(n) +2r(n)+21\delta+1$ satisfies the statement of the
lemma. \qed

\begin{cor}\label{sym}
Let $H$ be a torsion-free $\delta$-hyperbolic group with respect to a finite generating set $S$. There exists a
computable function $\hbar:\mathbb{N}\rightarrow \mathbb{R}^{+}$ with the following property: if $b,x_1,x_2,x_3 \in H$
and $0\neq m_1,m_2,m_3\in \mathbb{Z}$ are such that $||b^{m_1}||, ||b^{m_2}||, ||b^{m_3}||>15\delta$, $x_1x_2x_3=1$,
$m_1+m_2+m_3=0$, and $x_1b^{m_1}x_2b^{m_2}x_3b^{m_3}=1$, then each of the $x_i$ can be written in the form
$b^{n_1}ub^{n_2}vb^{n_3}$, where $n_1,n_2,n_3\in \mathbb{Z}$, and both $u,v$ have length at most $\hbar(|b|)$.
\end{cor}

\demo Inverting the last equation and cyclically permuting if necessary, we may assume that $m_1>0$ and $m_2>0$. Now,
Lemma~\ref{5} gives the conclusion. \qed

\medskip

We can now prove the following special case of Theorem~\ref{maingenbeg}.

\begin{prop}\label{cyclic}
Let $H$ be a torsion-free $\delta$-hyperbolic group with respect to a finite generating set $S$. Then, for any $g\in H$
there is a computable constant $C=C(|g|)>0$ with the following property: for every $a,b\in \langle g\rangle$ with
$||a||, ||b||, ||ab^{\pm 1}||>15\delta$, and every conjugate $b_*$ of $b$, if $ab_*^s$ is conjugate to $ab^s$ for every
$s=-C,\ldots ,C$, then $b_*=b$.
\end{prop}

\demo Let $a=g^n$ and $b=g^m$ (with $n,m\neq 0$ and $n\neq \pm m$), and let $b_{\ast}=x^{-1}bx$ for some $x\in H$
(which can always be multiplied on the left by a power of $b$).

We may assume $n,m>0$. Indeed, if $n<0$, we replace $g$ by $g^{-1}$, and $n$ by $-n$, and $m$ by $-m$; the statement
does not change and we get $n>0$. If then $m<0$, we replace $b$ by $b^{-1}=g^{-m}$ and $b_{\ast}$ by $b_{\ast}^{-1}$;
again the statement does not change and we get $m>0$.

So, let us assume $n,m>0$, $||a||, ||b||, ||ab^{\pm 1}||>15\delta$, and $ab_*^s$ being conjugate to $ab^s$ for every
$s=-C,\ldots ,C$, where $C$ is yet to be determined.

Taking $C\geqslant 1$, we have $ab_{\ast}^{-1}$ conjugate to $ab^{-1}$, that is $g^n\cdot x^{-1}g^{-m}x=h^{-1}g^{n-m}h$
for some $h\in H$. Rewrite this last equation into the following two forms
 \begin{equation}\label{first}
xh^{-1}g^{m-n}hx^{-1}\cdot xg^nx^{-1}=g^m,
 \end{equation}
 \begin{equation}\label{second}
h^{-1}g^{n-m}h\cdot x^{-1}g^mx=g^n.
 \end{equation}
If $m>n$, then from equation~(\ref{first}) and Lemma~\ref{5} we get
 $$
x=g^{p}vg^q
 $$
for some $p,q\in \mathbb{Z}$ and $v\in H$ with $|v|\leqslant \hbar(|g|)$. Otherwise, $m<n$ and then from
equation~(\ref{second}) and Lemma~\ref{5} we get the same expression for $x$. Replacing $x$ by $g^{-p}x$, we can assume
$p=0$, i.e. $x=vg^q$. And now, replacing $b_{\ast}$ by $g^qb_{\ast}g^{-q}$, which does not affect neither the
hypothesis nor the conclusion of the proposition (recall that both $a$ and $b$ are powers of $g$), we may assume that
$x=v$, $|v|\leqslant \hbar(|g|)$.

Let us impose that, $ab_{\ast}^s$ and $ab^s=g^{n+sm}$ are conjugate, for some positive value of $s$. By
Lemma~\ref{easy1}, there exists $z_s\in H$ such that
 \begin{equation}\label{third}
g^n\cdot x^{-1}g^{sm}x=ab_{\ast}^s=z_s^{-1}\underset{{}^c}{\cdot}g^{n+sm}\underset{{}^c}{\cdot}z_s,
 \end{equation}
where the constant $c$ depends only on $|g|$, $\delta$ and $\sharp S$. By Proposition~\ref{quasigeod2} and
Lemma~\ref{computable}, we can compute a constant $C_0$ such that $|g^{n+sm}|>2c+\delta$, for every $s\geqslant C_0$.
Taking at least this value for $C$, and using Lemma~\ref{cc} and Corollary~\ref{powers}, we deduce that
 $$
|g^n|+|g^{sm}|+2|x|\geqslant |ab_*^s|> |g^{n+sm}|+2|z_s|-(4c+2\delta)\geqslant |g^n|+|g^{sm}|-2\mu+2|z_s|-(4c+2\delta),
 $$
where $\mu=\mu(|g|)$ is the computable function from Corollary~\ref{useful}. Hence, $|z_s|\leqslant
\hbar(|g|)+\mu(|g|)+2c+\delta$.

Finally, take $C=C_0 +\sharp \mathcal{B}\big(\hbar(|g|)+\mu(|g|)+2c+\delta \big)$. Having $ab_*^s$ conjugate to $ab^s$
for every $s=-C,\ldots ,C$, we obtain elements $z_s$, $s=C_0,\ldots ,C$, all of them in the ball
$\mathcal{B}\Big(\hbar(|g|)+\mu(|g|)+2c+\delta \Big)$ by the previous paragraph.

Hence, there must be a repetition,
i.e. there exist $C_0<s_1<s_2<C$ such that $z_{s_1}=z_{s_2}$ (denote it by $z$). We have
 \begin{equation}
ab_{\ast}^{s_1}=z^{-1}g^{n+s_1m}z\label{vigo}
 \end{equation}
and
 $$
ab_{\ast}^{s_2}=z^{-1}g^{n+s_2m}z,
 $$
from which we deduce
 $$
b_{\ast}^{s_2-s_1}=z^{-1}g^{m(s_2-s_1)}z.
 $$
This implies $b_{\ast}=z^{-1}g^{m}z$, and then~(\ref{vigo}) implies $a=z^{-1}g^nz$. Since $a=g^n$, the element $z$
commutes with $g$ and so, again from~(\ref{vigo}), $b_{\ast}=b$. $\Box$

\section{The main theorem for two words}\label{2words}

The following lemma is a preliminary step in proving the main result for the case of two words (Theorem~\ref{main}).
Note that equations~(\ref{eq2}) and~(\ref{eq1}) in its formulation have the following common form: the product of
certain conjugates of two elements equals the product of these two elements.

\begin{lem} \label{twoeq}
Let $H$ be a torsion-free $\delta$-hyperbolic group with respect to a finite generating set $S$, and let $b,w\in H$.
There exists a computable constant $M=M(|b|,|w|)$ such that the following holds: if $b_*$ is conjugate to $b$ (say
$b_*=h^{-1}bh$), and $wb_*^k$ is conjugate to $wb^k$ for every $k=1,\ldots ,M$, then there exists an element $d\in H$
and integers $m,s,t$, such that $s+t>0$ and
 \begin{equation}\label{eq2}
(d\cdot b^s \cdot d^{-1})(dw \cdot b^t \cdot w^{-1}d^{-1})=b^{s+t},
 \end{equation}
 \begin{equation}\label{eq1}
(d^{-1}h\cdot w\cdot h^{-1}d)(d^{-1}\cdot b^m \cdot d) = wb^m.
 \end{equation}
\end{lem}

\demo The result is obvious if $b=1$. Let us assume $b\neq 1$.

If we prove the statement for a particular conjugator $h$, then we immediately have the same result for an arbitrary
other, just replacing $h$ to $b^qh$ and $d$ to $b^qd$ (for $q$ rational). So, we can choose our favorite $h$.

By~Lemma~\ref{easy1}, there exists a conjugator $h\in H$ such that, for any integer
$k\geqslant 0$, we have $b_*^k =h^{-1}\underset{^{c}}{\cdot}b^k\underset{^{c}}{\cdot}h$, where $c=\delta+\mu(|b|)$. Let us show
the result for this particular $h$. Since this expression remains valid while enlarging the constant $c$, we shall
consider it with $c=3\delta+\mu(|b|)+|w|+1$ in order to match with other calculations below. Thus,
 \begin{equation}\label{eq3}
wb_*^k=w(h^{-1}\underset{^{c}}{\cdot}b^k\underset{^{c}}{\cdot}h),
 \end{equation}
for every $k\geqslant 0$. Suppose that $wb_*^k$ is conjugate to $wb^k$ for every $k=1,\ldots ,M$, where $M$ is still to
be determined. Then, by Lemma~\ref{circ}, for each of these $k$'s there exist an element $e_k\in H$ and an integer
$l_k$, such that $0\leqslant l_k\leqslant k$ and
 \begin{equation}\label{eq4}
wb_*^k =e_k^{-1}\underset{^c}{\cdot}(b^{k-l_k}wb^{l_k})\underset{^c}{\cdot} e_k.
 \end{equation}

By Corollary~\ref{length}, and Proposition~\ref{quasigeod2} and Lemma~\ref{computable}, there exists a computable
constant $k_0 =k_0 (|b|,|w|)>0$ such that both $|b^{k-l_k}wb^{l_k}|$ and $|b^k|$ are bigger than $2c+\delta$ for all
$k\geqslant k_0$.

We introduce the following notation: for two sequences of elements $u_k\in H$ and $v_k\in H$ (where $k$ runs through a
subset of $\mathbb{N}$) we write $u_k\approx v_k$ if $|u_k^{-1}v_k|$ is bounded from above by a computable function,
depending on $\delta$, $\sharp S$, $w$, and $b$ only (so, in particular, not depending on $k$). The function will be
clear from the context. Similarly, we write $|u_k|\approx |v_k|$ if $||u_k|-|v_k||$ is bounded from above by a
computable function, depending on the same arguments.

Take $k\geqslant k_0$. Then from~(\ref{eq3}) and~(\ref{eq4}), and with the help of Lemma~\ref{cc}, we deduce
 $$
|wb^k_{\ast}|\approx 2|h|+|b^k|
 $$
and
 $$
|wb^{k}_{\ast}|\approx 2|e_k|+|b^{k-l_k}wb^{l_k}|\approx 2|e_k|+|b^k|,
 $$
where the last approximation is due to Corollaries~\ref{powers} and~\ref{length}. Therefore $|e_k|\approx |h|$.

Now we will prove that $e_k\approx h$. For that, we realize the right hand side of~(\ref{eq3}) in the Cayley graph
$\Gamma(H,S)$ as the path starting at 1 and consisting of 4 consecutive geodesics with labels equal in $H$ to the
elements $w$, $h^{-1}$, $b^k$, and $h$. Analogously, we realize the right hand side of~(\ref{eq4}) as the path starting
at 1 and consisting of 3 consecutive geodesics with labels equal in $H$ to the elements $e_k^{-1}$,
$b^{k-l_k}wb^{l_k}$, and $e_k$ (see Figure 3).

\unitlength 1mm % = 2.85pt
\linethickness{0.4pt}
\ifx\plotpoint\undefined\newsavebox{\plotpoint}\fi % GNUPLOT compatibility
\begin{picture}(100,26)(-5,20)
\qbezier(40,30)(45.25,31.25)(48.5,34.5) \qbezier(48.5,34.5)(58.38,32)(66.75,35.5)
\qbezier(66.75,35.5)(79,31.38)(90.25,33.75) \qbezier(90.25,33.75)(99,29.13)(109.75,30)
\qbezier(110,29.75)(98.38,27.63)(93.25,24) \qbezier(93.25,24)(70.38,28.88)(54,25.25)
\qbezier(54,25.25)(47.5,28.38)(40,30) \put(40,30){\circle*{1.5}} \put(48.5,34.25){\circle*{1.5}}
\put(67,35.25){\circle*{1.5}} \put(90.5,33.75){\circle*{1.5}} \put(110.25,30){\circle*{1.5}}
\put(93.5,23.75){\circle*{1.5}} \put(53.75,25){\circle*{1.5}} \put(43,33.5){\makebox(0,0)[cc]{$w$}}
\put(57,36){\makebox(0,0)[cc]{$h$}} \put(79.5,36){\makebox(0,0)[cc]{$b^k$}} \put(100.25,33){\makebox(0,0)[cc]{$h$}}
\put(46,26){\makebox(0,0)[cc]{$e_k$}} \put(74,23.5){\makebox(0,0)[cc]{$b^{k-l_k}wb^{l_k}$}}
\put(101.25,24.5){\makebox(0,0)[cc]{$e_k$}}
%\qbezier(93.5,24)(91.63,28)(93.25,32)
\put(36.5,30.5){\makebox(0,0)[cc]{$1$}} \put(118.75,30){\makebox(0,0)[cc]{$wb_{\ast}^k=${\small{$C$}}}}
\put(43.75,31.05){\vector(2,1){2.5}} \put(58,33.55){\vector(-1,0){0.5}} \put(80,33){\vector(1,0){0.5}}
\put(99,30.55){\vector(3,-1){0.5}} \put(47.5,28){\vector(-2,1){0.5}} \put(75,26.5){\vector(1,0){0.5}}
\put(100,27.25){\vector(2,1){0.5}} \put(40,30){\line(1,0){70}} \put(92,37){\makebox(0,0)[cc]{\small{$X$}}}
\put(94,20.5){\makebox(0,0)[cc]{\small{$Y$}}} \put(90.5,34){\line(0,-1){4.25}} \put(93.5,24){\line(0,1){6}}
\put(93.5,29.75){\circle*{1.5}} \put(90.5,29.75){\circle*{1.5}} \put(90,27.25){\makebox(0,0)[cc]{\small{$A$}}}
\put(95.5,27.5){\makebox(0,0)[cc]{\small{$B$}}}
\end{picture}

\begin{center} Figure 3 \end{center}

Both paths are $(\lambda,\epsilon)$-quasigeodesics connecting $1$ and $C=wb_{\ast}^k$, where $\lambda$ and $\epsilon$
are computable and depend only on $c$. We choose a geodesic $[1,C]$ and denote $X=wh^{-1}b^k$,
$Y=e_k^{-1}b^{k-l_k}wb^{l_k}$.

By Proposition~\ref{prop}, these quasigeodesics are both at bounded distance $R=R(\delta, c)$ from the segment $[1,C]$.
Therefore there are points $A,B\in [1,C]$, such that $|XA|\leqslant R$ and $|YB|\leqslant R$. In our notations we can
write $|XA|\approx 0$ and $|YB|\approx 0$. Therefore $|AC|\approx |XC|=|h|$ and $|BC|\approx |YC|=|e_k|$. Since
$|h|\approx |e_k|$, we have $|AC|\approx |BC|$ and so $|AB|\approx 0$. Hence, $|he_k^{-1}|=|XY|\leqslant
|XA|+|AB|+|BY|\approx 0$. This means that $e_k\approx h$ and so, $e_k$ lies in the ball with center $h$ and radius
depending only on $|b|$ and $|w|$.

Let $M$ be $1+k_0$ plus the number of elements in this ball. There must exist $k_0\leqslant k_1 <k_2 \leqslant M$ such
that $e_{k_1}=e_{k_2}$. Denote this element by $e$ and, rewriting equation~(\ref{eq4}) for these two special values of
$k$,
 \begin{equation}\label{averlange}
wb_*^{k_1} =e^{-1}(b^{k_1-l_{k_1}}wb^{l_{k_1}})e
 \end{equation}
and
 $$
wb_*^{k_2} =e^{-1}(b^{k_2-l_{k_2}}wb^{l_{k_2}})e,
 $$
we get
 $$
b_*^{k_2-k_1}= e^{-1}(b^{-l_{k_1}}w^{-1}b^{k_2-k_1+l_{k_1}-l_{k_2}}wb^{l_{k_2}})e.
 $$
Let $s=k_2-k_1+l_{k_1}-l_{k_2}$ and $t=l_{k_2}-l_{k_1}$ (so $s+t>0$). Recalling that $b_*^{k_2-k_1} =h^{-1}b^{k_2-k_1
}h$, we can rewrite the previous equation as
 $$
he^{-1} b^{-l_{k_1}}w^{-1}b^s w b^t b^{l_{k_1}} eh^{-1} = b^{s+t}.
 $$
Setting $d=he^{-1}b^{-l_{k_1}}w^{-1}$, we deduce $(db^s d^{-1})\cdot (dwb^t w^{-1}d^{-1}) =b^{s+t}$, which is
equation~(\ref{eq2}). And using equation~(\ref{averlange}), the definition of $d$ and $b_*^{k_1} =h^{-1}b^{k_1} h$, we
obtain $(d^{-1}hwh^{-1}d)\cdot (d^{-1}b^{k_1} d) =wb^{k_1}$, which is equation~(\ref{eq1}) with $m=k_1$. \qed

Now, using~(\ref{eq2}) and~(\ref{eq1}) and distinguishing the cases $st\neq 0$ or $st=0$, we will obtain more
information about relations between $w,b$ and $h$.

\begin{prop}\label{non-zero}
Let $H$ be a torsion-free $\delta$-hyperbolic group with respect to a finite generating set $S$ and let $b,w,d$ be
elements of $H$ satisfying equation~{\rm(\ref{eq2})}. Suppose additionally that $||b^k ||>15\delta$ for all $k>0$, and
that $st\neq 0$. Then, there exist integers $p,q,r$ and elements $u,v\in H$ of length at most $\hbar(|b|)$, such that
 $$
w=b^pub^rvb^q.
 $$
\end{prop}

\demo This follows directly from Corollary~\ref{sym}. \qed

\begin{prop}\label{zero}
Let $H$ be a torsion-free $\delta$-hyperbolic group with respect to a finite generating set $S$ and let $b,w,d,h$ be
elements of $H$ satisfying equations~{\rm(\ref{eq2})} and~{\rm(\ref{eq1})} with $s+t>0$. Suppose additionally that
$st=0$. Then $h=b^pw^q$ for some rational numbers $p,q$.
\end{prop}

\demo Let us distinguish two cases.

{\it Case 1: $s=0$}. In this case, equation~(\ref{eq2}) says that $dw$ commutes with $b$. So, $dw=b^p$ for some
rational $p$. Plugging this into equation~(\ref{eq1}) we obtain $hwh^{-1}=b^{p+m}wb^{-p-m}$. Hence, $b^{-p-m}h$
commutes with $w$ and the result follows.

{\it Case 2: $t=0$}. In this case, equation~(\ref{eq2}) says that $d$ commutes with $b$. So, $d=b^p$ for some rational
$p$. Plugging this into equation~(\ref{eq1}) we obtain $b^{-p}hwh^{-1}b^p=w$. Hence, $b^{-p}h$ commutes with $w$ and
the result follows. $\Box$

\medskip

Next, we need to obtain some extra information by applying Lemma~\ref{twoeq} to sufficiently many different elements
$w$. To achieve this goal, given a pair of elements $a,b\in H$, we consider the finite set
 $$
\mathcal{W}=\{(a^i b)^{2j}\,|\, 1\leqslant i\leqslant 1+N,\,\, 1\leqslant j\leqslant 1+3N^2 \} \subseteq \langle
a,b\rangle \leqslant H,
 $$
where
 $$
N=N(|b|)=\sharp \mathcal{B}\, (\hbar(|b|)),
 $$
and $\hbar$ is the function from Lemma~\ref{5}. Let us systematically apply Lemma~\ref{twoeq} to every $w\in
\mathcal{W}$.

\begin{lem}\label{st}
Let $H$ be a torsion-free $\delta$-hyperbolic group with respect to a finite generating set $S$. Let $a,b\in H$ be
elements generating a free subgroup of rank 2, and with $||b^k ||>15\delta$ for all $k>0$. Suppose that for every $w\in
\mathcal{W}$, there exists a conjugate $b_*$ of $b$ such that the elements $w,b,b_*$ satisfy the hypothesis of
Lemma~\ref{twoeq} (i.e. $wb_*^k$ is conjugate to $wb^k$, for every integer $k=1,\dots, M(|b|,|w|)$). Then, for at least
one such $w\in \mathcal{W}$, the conclusion of Lemma~\ref{twoeq} holds with $st=0$.
\end{lem}

\demo Under the hypothesis of the lemma, suppose that we have equations~{\rm(\ref{eq2})} and~{\rm(\ref{eq1})} with
$st\neq 0$ for every $w\in \mathcal{W}$, and let us find a contradiction.

Write $\mathcal{W}={\bigsqcup}_{i=1}^{1+N} \mathcal{W}_i$, where $\mathcal{W}_i=\{(a^i b)^{2j}\,|\, 1\leqslant
j\leqslant 1+3N^2 \}$, and fix a value for $i\in \{1,\ldots ,N+1\}$.

By Proposition~\ref{non-zero}, for every $w\in \mathcal{W}_i$, there exist integers $p,q,r$, and elements $u,v\in H$ of
length at most $\hbar(|b|)$ such that
 \begin{equation}\label{ijk}
b^pwb^q=ub^rv.
 \end{equation}
(of course, these integers and elements depend on $w$). Since $\sharp {\mathcal W}_i =1+3N^2>3(\sharp
\mathcal{B}(\hbar(|b|)))^2$ (because $\langle a,b\rangle$ is free of rank 2) and the lengths of $u$ and $v$ are at most
$\hbar(|b|)$, there must exist four diferent elements of $\mathcal{W}_i$ with the same $u$ and $v$. That is, there
exists $w_1=(a^i b)^{\sigma}$, $w_2=(a^i b)^{\tau}$, $w_3=(a^i b)^{\sigma'}$ and $w_4=(a^i b)^{\tau'}$ (where the
exponents $0<\sigma<\tau<\sigma'<\tau'$ all differ at least 2 from each other) such that
 $$
\begin{array}{ccc} b^{p_1}w_1 b^{q_1}=ub^{r_1}v, & \quad & b^{p_2}w_2 b^{q_2}=ub^{r_2}v, \\ b^{p_3}w_3
b^{q_3}=ub^{r_3}v, & \quad & b^{p_4}w_4 b^{q_4}=ub^{r_4}v. \end{array}
 $$
Combining these equations, we get
 \begin{equation}\label{ik}
b^{p_2}w_2b^{q_2-q_1}w_1^{-1}b^{-p_1}=ub^{r_2-r_1}u^{-1},
 \end{equation}
 $$
b^{p_4}w_4b^{q_4-q_3}w_3^{-1}b^{-p_3}=ub^{r_4-r_3}u^{-1}.
 $$
Hence, the left hand sides of these two equations commute. Let us rewrite them in the form
 $$
x=b^{\alpha}(a^i b)^{\tau} b^{\beta}(a^i b)^{-\sigma} b^{\gamma},
 $$
 $$
x'=b^{\alpha'}(a^i b)^{\tau'} b^{\beta'}(a^i b)^{-\sigma'} b^{\gamma'},
 $$
where $0<\sigma<\tau$ and $0<\sigma'<\tau'$ all differ at least 2 from each other (and we have no specific information
about the integers $\alpha, \beta, \gamma, \alpha', \beta', \gamma'$). The key point here is that this commutativity
relation between $x$ and $x'$ happens inside the free group $\langle a,b\rangle$.

Consider now the monomorphism $\langle a,b\rangle \to \langle a,b\rangle$ given by $a\mapsto a^i b$, $b\mapsto
b$. Since $x$ and $x'$ both lie in its image, and commute, their preimages, namely $y=b^{\alpha} a^{\tau}
b^{\beta} a^{-\sigma} b^{\gamma}$ and $y'=b^{\alpha'} a^{\tau'} b^{\beta'} a^{-\sigma'} b^{\gamma'}$,
must also commute.

Suppose $\beta \beta' \neq 0$. Then, $y$ is not a proper power in $\langle a,b\rangle$ (in fact, its cyclic
reduction is either $a^{\tau} b^{\beta} a^{-\sigma} b^{\alpha +\gamma}$ with $\alpha+\gamma \neq 0$, or
$a^{\tau-\sigma} b^{\beta}$, which are clearly not proper powers). Similarly, $y'$ is not a proper power either.
Then the commutativity of $y$ and $y'$ forces $y=y'^{\pm 1}$, which is obviously not the case. Hence, $\beta \beta'
=0$. Without loss of generality, we can assume $\beta =0$.

Let us go back to equation~(\ref{ik}) which, particularized to this special case, reads
 $$
b^{\alpha }(a^i b)^{\tau} b^0 (a^i b)^{-\sigma} b^{\gamma }=u b^{\delta} u^{-1},
 $$
that is
\begin{equation}\label{k}
b^{\alpha }(a^i b)^{\rho} b^{\gamma }=ub^{\delta} u^{-1},
 \end{equation}
where $\rho=\tau-\sigma \geqslant 2$.
Recall that all these arguments were started for a fixed value of $i$ and that
the corresponding element $u$ (which depends on the chosen $i$) has length at most $\hbar(|b|)$.

Finally, it is time to move $i=1,\ldots ,1+N$. Since $1+N>\sharp \mathcal{B}(\hbar(|b|))$, there must exist two indices
$1\leqslant i_1<i_2\leqslant 1+N$ giving the same $u$. Equation~(\ref{k}) in these two special cases is
 $$
b^{\alpha }(a^{i_1} b)^{\rho} b^{\gamma } = u b^{\delta} u^{-1}
 $$
and
 $$
b^{\alpha' }(a^{i_2} b)^{\rho'} b^{\gamma' } = u b^{\delta'} u^{-1},
 $$
where $\rho, \rho' \geqslant 2$ and $1\leqslant i_1<i_2$. Again, $z=b^{\alpha }(a^{i_1} b)^{\rho} b^{\gamma }$ and
$z'=b^{\alpha' }(a^{i_2} b)^{\rho'} b^{\gamma'}$ commute. Since $i_1, i_2, \rho$ and $\rho'$ are all positive, this
implies that some positive power of $z$ equals some positive power of $z'$. But it is straightforward to see that
(after all possible reductions) the first $a$-syllable of any positive power of $z$ is $a^{i_1}$ (here we use
$\rho\geqslant 2$); similarly the first $a$-syllable of any positive power of $z'$ is $a^{i_2}$. Since $i_1 \neq i_2$,
this is a contradiction and the proof is completed. \qed

\medskip

Now can already prove the main Theorem~\ref{maingenbeg}, in the special case $n=2$.

\begin{thm}\label{main}
Let $H$ be a torsion-free $\delta$-hyperbolic group with respect to a finite generating set $S$, and consider four
elements $a,b, a_*, b_*\in H$ such that $a_*$ is conjugate to $a$, and $b_*$ is conjugate to $b$. There exists a
computable constant $L$ (only depending on $|a|$, $|b|$, $\delta$ and $\sharp S$), such that if $(a_*^i b_*^l)^{j}
b_*^k$ is also conjugate to $(a^i b^l)^{j} b^k$ for every $i,j,k,l=-L,\ldots ,L$ then there exists a uniform conjugator
$g\in H$ with $a_*=g^{-1}ag$ and $b_*=g^{-1}bg$ (i.e. $(a_*,b_*)$ is conjugate to $(a,b)$).
\end{thm}

\demo The conclusion is obvious if  $a$ or $b$ is trivial. So, let us assume $a\neq 1$ and $b\neq 1$.
Note, that $\langle a\rangle=\langle b\rangle$  and even $a=b^{\pm 1}$ is allowed.

Suppose that $(a_*^i b_*^l)^{j} b_*^k$ is conjugate to $(a^i b^l)^{j} b^k$ for every $i,j,k,l=-L,\ldots ,L$, where $L$
is still to be determined. We shall prove the result imposing several times that $L$ is big enough, in a constructive
way. At the end, collecting together all these requirements, we shall propose a valid value for $L$.

Since $H$ is torsion-free, every nontrivial element has infinite cyclic centralizer (see Proposition~\ref{nova}). Let
$a_1,b_1$ be generators of $C_H(a)$ and $C_H(b)$. Inverting $a_1$ or $a_2$ if necessary, we may assume that $a=a_1^p$
and $b=b_1^q$ for positive $p$ and $q$. By Corollary~\ref{2}, there exists a computable natural number $r_0$ such that
for every $r\geqslant r_0$, $||a_1^r||>15\delta$ and $||b_1^r||>15\delta$. So, after replacing $a,b,a_*,b_*$ by
$a^{r_0},b^{r_0},a_*^{r_0},b_*^{r_0}$, we can assume that $||a^r||>15\delta$ and $||b^r||>15\delta$ for every $r\neq
0$. Moreover, if $a,b$ generate a cyclic group, then after the above replacement either $a=b$ or
$||ab^{-1}||>15\delta$. Analogously, either $a=b^{-1}$, or $||ab||>15\delta$.

For every word $w$ on $a$ and $b$, let us denote by $w_*$ the corresponding word on $a_*$ and $b_*$. Now, observe that
we can uniformly conjugate $a_*$ and $b_*$ by any element of $H$ (and abuse notation denoting the result $a_*$ and
$b_*$ again), and both the hypothesis and conclusion of the theorem does not change. In particular, for any chosen word
of the form $w=(a^i b^l)^{j} b^k$ (with $i,j,k,l=-L,\ldots ,L$), we can assume that $w_*=w$ (of course, with an
underlying $a_*$ and $b_*$ now depending on $w$); when doing this, we say that we \emph{center the notation on $w$}.
Note that centering notation does not change $a,b$, therefore the constant $L$ is not affected.

Let us distinguish two cases.

\emph{Case 1: $\langle a,b\rangle$ is a cyclic group, say $\langle g\rangle$}. Centering the notation on $a$, we may
assume that $a_*=a$. If $a=b^{\epsilon}$, where $
\epsilon =\pm 1$, then we use that $ab_{\ast}^{-\epsilon}$ is conjugate to $ab^{-\epsilon}=1$ and deduce immediately
that $b_{\ast}=b$. Now, assume that $a\neq b^{\pm 1}$, and so $||ab^{\pm 1}||>15\delta$. Part of our hypothesis says that
$a_*b_*^l =ab_*^l$ is conjugate to $ab^l$ for every $l=-L,\ldots ,L$. Hence, taking $L$ bigger than or equal to the
constant $C=C(|g|)$ from Proposition~\ref{cyclic}, we obtain $b_*=b$. This concludes the proof in this case.

\emph{Case 2: $\langle a,b\rangle$ is not cyclic}. By Proposition~\ref{free}, there exists a sufficiently big and
computable natural number $p$ such that $\langle a^p,b^p\rangle$ is a free subgroup of $H$ of rank 2. Note that,
multiplying the constant by $p$, and using the uniqueness of root extraction in $H$, the result follows from the same
result applied to the elements $a^p, b^p$ and $a_*^p, b_*^p$. So, after replacing $a,b,a_*,b_*$ by
$a^p,b^p,a_*^p,b_*^p$, we can assume that $F_2 \simeq \langle a,b\rangle \leqslant H$.

With these gained assumptions, let us show that any constant
 $$
L\geqslant \max\{2+6N^2, \underset{w\in \mathcal{W}}{\max}\, M(|b|,|w|)\},
 $$
works for our purposes, where the number $N$ and the set $\mathcal{W}$ are defined before Lemma~\ref{st}, and the
function $M$ is defined in Lemma~\ref{twoeq}.

Part of our hypothesis says that, for every $w=(a^ib)^{2j}\in \mathcal{W}$, $w_*b_*^k=(a_*^ib_*)^{2j}b_*^k$ is
conjugate to $wb^k$ for every $k=1,\ldots ,M(|b|,|w|)$.

Fix $w\in \mathcal{W}$. Centering the notation on this $w$, we have that $wb_*^k$ $(=w_*b_*^k)$ is conjugate to $wb^k$
for every $k=1,\ldots ,M(|b|,|w|)$. That is, $w$ satisfies the hypothesis of Lemma~\ref{twoeq} (with the corresponding
value of $b_*$). And this happens for every $w\in \mathcal{W}$. Thus, Lemma~\ref{st} ensures us that the conclusion of
Lemma~\ref{twoeq} holds with $st=0$ for at least one $w_0 =(a^{i_0} b)^{2j_0}\in \mathcal{W}$, $1\leqslant i_0\leqslant
1+N$, $1\leqslant j_0\leqslant 1+3N^2$ (note that Lemma~\ref{st} can be applied because we previously gained the
assumptions $||b^r||>15\delta$ for every $r\neq 0$, and $F_2 \simeq \langle a,b\rangle \leqslant H$). For the rest of
the proof, let us center the notation on this particular $w_0$.

Using Proposition~\ref{zero}, we conclude that every conjugator from $b$ to $b_*$ (say $b_*=h^{-1}bh$) is of the form
$h=b^pw_0^q$ for some rational numbers $p,q$. Hence, $w_0^{-q}bw_0^q =b_*$. Then,
 $$
((w_0^{-q}aw_0^q)^{i_0}b_*)^{2j_0}=w_0^{-q}(a^{i_0}b)^{2j_0} w_0^q =w_0^{-q}w_0 w_0^q =w_0 =w_{0*} =(a_*^{i_0}b_*)^{2j_0}.
 $$
Extracting roots twice, we conclude that $w_0^{-q}aw_0^q =a_*$. Thus, $w_0^q$ is a uniform right conjugator from
$(a,b)$ to $(a_*, b_*)$. This concludes the proof for this second case. \qed

\section{Main theorem for several words}

Finally, we extend the result to arbitrary tuples of words, thus proving the main result of the paper.

\emph{Proof of Theorem~\ref{maingenbeg}.} The implication to the right is obvious (without any bound on the length of
$W$).

Let $\mathcal{A}=\{a_1,\ldots,a_n\}$, and assume that $W(a_{1*},\ldots ,a_{n*})$ is conjugate to $W(a_1,\ldots ,a_n)$
for every word $W$ in $n$ variables and length up to a constant yet to be determined. As above, we shall prove the
result assuming several times this constant to be big enough, in a constructive way. The reader can collect together
all these requirements, and find out a valid explicit value (which will depend only on $\delta$, $\sharp{S}$ and
$\sum_{i=1}^n |a_i|$). Decreasing $n$ if necessary, we may assume that all $a_i$ are nontrivial. If $n=1$ there is
nothing to prove, so assume $n\geqslant 2$.

Suppose the elements $a_1,\ldots,a_n$ generate a cyclic group, say $\langle a_1,\ldots ,a_n\rangle \leqslant \langle
g\rangle \leqslant H$, with $g$ root-free. Applying Theorem~\ref{main} to every pair $a_1,a_j$, we get a computable
constant such that if $W(a_{1*}, a_{j*})$ is conjugate to $W(a_1,a_j)$ for every word $W$ of length up to this
constant, then $a_1$ and $a_j$ admit a common conjugator, say $x_j$. Taking the maximum of these constants over all
$j=2,\ldots ,n$ we are done, because $x_j^{-1}a_1x_j=a_{1*}$ and $x_j^{-1}a_jx_j=a_{j*}$ for $j=2,\ldots ,n$ imply that
$x_2x_j^{-1}\in C_H(a_1)=\langle g\rangle$, and hence $x_2^{-1}a_jx_2 = x_j^{-1}(x_jx_2^{-1}a_jx_2x_j^{-1}) x_j =
x_j^{-1} a_j x_j =a_{j*}$ for $j=2,\ldots ,n$; thus, $x_2$ becomes a common conjugator.

So, we are reduced to the case where two elements of $\mathcal{A}$, say $a_1$ and $a_2$, generate a noncyclic group. In
this case, by Proposition~\ref{free}, there is a big enough computable $m$ such that $\langle a_1^m,a_2^m\rangle$ is a
free group of rank 2. Replacing $a_1, a_2$ by $a_1^m, a_2^m$ and $a_{1*}, a_{2*}$ by $a_{1*}^m, a_{2*}^m$, and
multiplying the computable constant by $m$, we may assume that $\langle a_1, a_2\rangle$ is free of rank 2.

By Theorem~\ref{main} (and taking the constant appropriately big), $a_1$ and $a_2$ admit a common conjugator. So,
conjugating the whole tuple $a_{1*}, \ldots ,a_{n*}$ accordingly, we may assume that $a_{1\ast}=a_1$ and
$a_{2\ast}=a_2$. We will prove that $a_{j\ast}=a_j$ for every $j=3,\ldots n$ as well.

By Lemma~\ref{minasyan} twice, there exists a big enough computable $k\geqslant 2$ such that the elements $a_1a_2^k$
and $a_2(a_1a_2^k)^k$ are root-free (and form a new basis for $\langle a_1,a_2\rangle$). Replacing $a_1$ by $a_1a_2^k$
and $a_{1*}$ by $a_{1*}a_{2*}^k$, and $a_2$ by $a_2(a_1a_2^k)^k$ and $a_{2*}$ by $a_{2*}(a_{1*}a_{2*}^k)^k$, and
updating the constant, we may assume that both $a_1$ and $a_2$ are root-free in $H$.

For every $j\geqslant 3$, let us apply Theorem~\ref{main} to the pairs $(a_1,a_j)$ and $(a_{1*}=a_1, a_{j*})$; we
obtain $x_j\in C_H(a_1)=\langle a_1\rangle$ such that $a_{j\ast}=x_j^{-1}a_jx_j$. Analogously, playing with the pair of
indices $2,j$, we get $y_j\in C_H(a_2)=\langle a_2\rangle$ such that $a_{j\ast}=y_j^{-1}a_jy_j$. In particular,
$x_j=a_1^{p_j}$ and $y_j=a_2^{q_j}$ for some integers $p_j,q_j$. Furthermore, $x_jy_j^{-1}\in C_{H}(a_j)$, that is
$a_1^{p_j}a_2^{-q_j}=a_j^{r_j}$ for some rational $r_j$. Note that if $p_jq_j=0$ then $a_{j*}=a_j$ as we want.

Again by Lemma~\ref{minasyan}, there is a big enough computable $k'\geqslant 2$ such that $b_1=a_1a_2^{k'}$ and
$b_2=a_2(a_1a_2^{k'})^{k'}$ are again root-free in $H$. Arguing like in the previous paragraph with these new elements,
we deduce a similar conclusion: for each $j=3,\ldots ,n$, either $a_{j\ast}=a_j$, or $b_1^{p_j'}b_2^{-q_j'}=a_j^{r_j'}$
for some nonzero integers $p_j',q_j'$ and some rational $r_j'$.

Thus, for each $j=3,\ldots ,n$, we either have (1) $a_{j\ast}=a_j$, or (2) $a_1^{p_j}a_2^{-q_j}=a_j^{r_j}$ and
$b_1^{p_j'}b_2^{-q_j'}=a_j^{r_j'}$ for some nonzero integers $p_j, q_j, p_j', q_j'$ and some rationals $r_j, r_j'$. But
this last possibility would imply that the elements $a_1^{p_j}a_2^{-q_j}$ and $b_1^{p_j'}
b_2^{-q_j'}=(a_1a_2^{k'})^{p_j'}(a_2(a_1a_2^{k'})^{k'})^{-q_j'}$ commute in the free group $\langle a_1,a_2\rangle$,
which is not the case, taking into account that $p_jq_jp_j'q_j'k'\neq 0$. Therefore,
$a_{j\ast}=a_j$ for each $j=1,\ldots,n$ and the proof is complete. \qed

\section{A mixed version for Whitehead's algorithm}

Particularizing the main result of the paper to the case of finitely generated free groups, we will obtain a mixed
version of Whitehead's algorithm.

Let us consider lists of elements in a finitely generated free group $F$, organized in $n$ blocks:
 $$
u_{1,1},\ldots ,u_{1,m_1}\,\,;\,\, \ldots \,\,;\,\, u_{i,1},\ldots ,u_{i,m_i}\,\,;\,\,\ldots \,\,;\,\,u_{n,1}, \ldots
,u_{n,m_n}.
 $$
The \emph{mixed Whitehead problem} consists in finding an algorithm to decide whether, given two such lists, there
exists an automorphism of $F$ sending the first list to the second up to conjugation, but asking for a uniform
conjugator in every block (and possibly different from those in other blocks).

Note that in the case where each block consists of one element (i.e. $m_i=1$ for all $i=1,\ldots ,n$), this is exactly
asking whether there exists an automorphism of $F$ sending the first list of elements to the second one up to
conjugacy, with no restriction for the conjugators. This problem (we call it the {\it Whitehead problem} for $F$) was already solved by Whitehead back in 1936 (see~\cite{W}
or~\cite{LS}).

On the other hand, if there is only one block (i.e. $n=1$), the problem is equivalent to ask whether there exists an
automorphism of $F$ sending the first list of elements exactly to the second. This was solved in 1974 by McCool
(see~\cite{Mc} or~\cite{LS}).

As a corollary of Theorem~\ref{maingenbeg}, we deduce a solution to the mixed Whitehead problem.

\begin{thm}\label{Whitehead}
Let $F$ be a finitely generated free group. Given two lists of words in $F$, $u_{i,j}$ and $v_{i,j}$, for $i=1,\ldots
,n$ and $j=1,\ldots ,m_i$, it is algorithmically decidable whether there exists $\varphi \in {\text{\rm Aut}}(F)$ and
elements $z_i \in F$ such that $\varphi(u_{i,j})=z_i^{-1}v_{i,j}z_i$ for every $i=1,\ldots ,n$ and $j=1,\ldots ,m_i$.
\end{thm}

\demo For every $i=1,\ldots ,n$, we compute the constant $C_i$ (depending only on $\sum_{j=1}^{m_i} |u_{i,j}|$
and the ambient rank) given in Theorem~\ref{maingenbeg} for the tuples of words $u_{i,1},\ldots ,u_{i,m_i}$ and
$v_{i,1},\ldots ,v_{i,m_i}$. By Theorem~\ref{maingenbeg}, an automorphism $\alpha\in {\text{\rm Aut}}(F)$ sends each
$W(u_{i,1},\ldots ,u_{i,m_i})$ to a conjugate of $W(v_{i,1},\ldots ,v_{i,m_i})$ (for every $W$ of length less than or
equal to $C_i$), if and only if $\alpha$ sends each $u_{i,j}$ to $z_i^{-1}v_{i,j}z_i$, $j=1,\ldots ,m_i$, for some
uniform conjugator $z_i$.

Now, let us enlarge each block of $u$'s and $v$'s with all the words of the form $W(u_{i,1},\ldots ,u_{i,m_i})$ and
$W(v_{i,1},\ldots ,v_{i,m_i})$, respectively, where $W$ runs over the set of all words in $m_i$ variables and length
less than or equal to $C_i$. Our problem is now equivalent to deciding whether there exists an automorphism $\varphi
\in {\text{\rm Aut}}(F)$ sending $W(u_{i,1},\ldots ,u_{i,m_i})$ to a conjugate of $W(v_{i,1},\ldots ,v_{i,m_i})$ for
every $i$, and for every $W$ of length less than or equal $C_i$. This is decidable by the classical version of
Whitehead's algorithm. \qed

This proof shows that the following theorem is true.

\begin{thm}\label{WhiteheadHyp} Let $H$ be a torsion-free hyperbolic group. If the Whitehead problem for $H$ is solvable, then the mixed Whitehead problem for $H$ is also solvable.
\end{thm}

\section{\bf Acknowledgements}

The first named author thanks the MPIM at Bonn for its support and excellent working conditions during the fall 2008,
while this research was finished. The second named author gratefully acknowledges partial support from the MEC (Spain)
and the EFRD (EC) through projects number MTM2006-13544 and MTM2008-01550.

\end{document}